\documentclass[11pt]{article}

\usepackage{makeidx}
\usepackage{sec}
\usepackage{pb-diagram}
\usepackage{amssymb}
\usepackage{amsmath}
\usepackage{amsthm}

\theoremstyle{plain}
\newtheorem{theorem}{Theorem}[section]
\newtheorem{lemma}[theorem]{Lemma}
\newtheorem{proposition}[theorem]{Proposition}
\newtheorem{corollary}[theorem]{Corollary}
\newtheorem{definition}[theorem]{Definition}

\theoremstyle{definition}
\newtheorem{ex}[theorem]{Example}
\newtheorem{remark}[theorem]{Remark}

\DeclareMathOperator{\End}{{\rm End}}

\DeclareMathOperator{\tr}{{\rm tr}}

\DeclareMathOperator{\Vect}{{\rm Vect}}

\newcommand{\td}{\tilde{d}}
\newcommand{\te}{\tilde{e}}
\newcommand{\tf}{\tilde{f}}
\newcommand{\tg}{\tilde{g}}
\newcommand{\st}{\sqrt{t}}

\newcommand{\uc}{\underline{\mathbb C}}

\newcommand{\bA}{\mathbb A}

\newcommand\bC{{\mathbb C}}

\newcommand\bR{{\mathbb R}}
\newcommand\bS{{\mathbb S}}

\newcommand\bZ{{\mathbb Z}}

\renewcommand{\Box}{\blacksquare}

\newcommand{\gA}{{\mathfrak A}}

\newcommand{\gb}{{\mathfrak b}}

\newcommand{\gH}{{\mathfrak H}}

\newcommand{\gM}{{\mathfrak M}}

\newcommand{\gT}{{\mathfrak T}}

\newcommand{\gw}{{\mathfrak w}}

\newcommand{\dir}{{\mathfrak D}}

\newcommand{\e}{\mathcal E}

\newcommand{\h}{\mathcal H}
\newcommand{\J}{\mathcal J}

\newcommand{\R}{\mathcal R}

\newcommand{\ra}{\rightarrow}

\newcommand{\Lra}{{\longrightarrow}}

\newcommand{\lan}{\langle}
\newcommand{\ran}{\rangle}

\def\inpr{\mathbin{\hbox to 6pt{\vrule height0.4pt width5pt depth0pt \kern-.4pt \vrule height6pt width0.4pt depth0pt\hss}}}

\newcommand{\hd}{\hat{d}}

\newcommand{\hg}{\hat{g}}
\newcommand{\bc}{{\bf c}}
\newcommand{\hbc}{\hat{\bf c}}

\newcommand{\hN}{{\hat{N}}}

\newcommand{\ii}{{\bf i}}

\newcommand{\si}{{\sigma}}

\newcommand{\ve}{{\varepsilon}}

\newcommand{\pa}{\partial}
\newcommand{\bpar}{\bar{\partial}}

\newcommand{\nai}{\nabla_i}
\newcommand{\naj}{\nabla_j}
\newcommand\nak{\nabla_k}
\newcommand{\nah}{\hat{\nabla}}
\newcommand{\tna}{\tilde{\nabla}}

\newcommand{\uso}{\underline{so}}

\begin{document}

\title{Geometric  connections and geometric Dirac operators on contact manifolds}

\author{Liviu I. Nicolaescu\thanks{This work was partially supported by NSF grant DMS-0071820.}\\ Dept.of Mathematics\\
 University of Notre Dame\\
Notre Dame, IN 46556; USA\\ nicolaescu.1@nd.edu}


\date{Version 2: January 16,2001}

\maketitle

\section*{Introduction}

Suppose $(M^{2n+1}, \eta, g, J)$ is a positively oriented,  metric contact manifold. More precisely, this means   that  $\eta$ is a $1$-form  such that $\frac{1}{n!}\eta \wedge (d\eta)^n=dv_g$, $J$ is a skew-symmetric endomorphism of $TM$ such that
\[
J^2X=-X+ \eta(X) \xi,\;\;d\eta(X,Y)=g(JX,Y),\;\;\forall X,Y\in \Vect(M),
\]
and $\xi$ is the Reeb vector field determined by $\eta(\xi)=1$, $\xi\inpr d\eta =0$. Set $V=\ker\eta$.

The operator $J$ induces  an almost complex structure on $V$, and we  get  decompositions
\[
V\otimes {\bC}=V^{1,0}\oplus V^{0,1},\;\;\Lambda^*V^*\otimes {\bC}=\bigoplus_{0\leq p+q\leq 2n} \Lambda^{p,q}V^*.
\]
We set $\Omega^{p,q}(V^*):=C^\infty(\Lambda^{p,q}V^*)$.  The Lie derivative along $\xi$ has the property $L_\xi \Omega^{0,p}(V^*)\subset \Omega^{0,p}(V^*)\oplus \Omega^{1,p-1}(V^*)$, and we define $L^V_\xi:\Omega^{0,p}(V^*)\ra \Omega^{0,p}(V^*)$  by $L^V_\xi\phi= (L_\xi\phi)^{0,p}$. The operator $\ii L^V_\xi$ is symmetric. There exists a natural operator
\[
\bpar_V:\Omega^{0,*}(V^*)\ra \Omega^{0,*+1}(V^*).
\]
We can form a {\em contact Hodge-Dolbeault} operator
\[
\h :  \Omega^{0,*}(V^*)\ra \Omega^{0,*}(V^*)
\]
which with respect to the decomposition $\Omega^{0,even}\oplus \Omega^{0,odd}(V^*)$ has the block form
\[
\h =\left[
\begin{array}{cc}
-\ii L^V_\xi & \sqrt{2}(\bpar_V+\bpar_V^*)\\
 & \\
\sqrt{2}(\bpar_V^*+\bpar_V^*) & \ii L^V_\xi
\end{array}
\right]
\]
This is a symmetric Dirac type operator and it is  an example of geometric Dirac operator, i.e.   an operator of Dirac type  defined entirely in geometric terms with no mention  of $spin^c$ structures.

On the other hand, the  contact form  defines a $spin^c$ structure   with determinant line $K_M^{-1}$, where canonical line bundle  of $M$ is defined by
\[
K_M:=\det V^{0,1}\cong \Lambda^{n,0}V^*.
\]
The associated bundle of complex spinors  is ${\bS}_c=\Lambda^{0,*}V^*$. The Clifford multiplication by $\ii \eta$ is  an involution of ${\bS}_c$ and the $\pm 1$ eigenspaces are
\[
{\bS}_c^\pm\cong \Lambda^{0,even/odd}V^*.
\]
A metric connection $\nabla$ on $TM$   such that $\nabla J=0$  is called a {\em contact connection}.  If   additionally $M$ is a $CR$ manifold, i.e. the distribution $V^{1,0}$ is integrable, then we define a  $CR$ connection to be a contact connection  such that its torsion satisfies
\[
T(X,Y)=0\;\;\;\forall X,Y\in C^\infty(V^{1,0}).
\]
A metric connection  on $TM$ together with a hermitian $A$ connection on $K_M^{-1}$ canonically define a Dirac operator $\dir(\nabla,A)$ on ${\bS}_c$.  The connection $\nabla$ is called {\em nice} if $\dir(\nabla, A)$ is symmetric  for any  hermitian connection $A$ on $K_M^{-1}$.  Two  metric connections $\nabla^1$  and $\nabla^2$ are called {\em Dirac equivalent} if there exists a hermitian connection $A$ on $K_M^{-1}$ such that $\dir(\nabla^1,A)=\dir(\nabla^2,A)$.

The first question we address in this paper is the following.

\medskip

\noindent $\bullet$ {\it Can we find a contact connection $\nabla$ and a hermitian connection $A$ on $K_M^{-1}$  $\dir(\nabla, A)=\h$?}  ( A connection $\nabla$  with this property is said to be {\em adapted to} $\h$.)

\medskip

Suppose additionally that $M$ is also $spin$.  We denote by $\dir_0$ the associated $spin$ Dirac operator.   The second question we as is the following.

\medskip

\noindent $\bullet$ {\it Does there  exist a contact connection $\nabla$ on $TM$ such that  $\dir(\nabla)=\dir_0$? In other words,  is there any contact connection Dirac equivalent to the Levi-Civita connection?}

\medskip

To address these questions we rely on the work P. Gauduchon, (see \cite{Gau} or \ref{ss: 21},\ref{ss: 22}), concerning    hermitian connections on almost-hermitian manifolds.    We can  naturally associate  two almost hermitian manifolds to $M$.

\medskip

\noindent $\bullet$ The cylinder $\hat{M}={\bR}\times M$ with metric $\hg=dt^2+g$ and almost complex structure  $\hat{J}$ defined  by $\hat{J}\pa_t=\xi$, $\hat{J}\mid_V=J$.

\medskip

\noindent $\bullet$ The symplectization $\tilde{M}= {\bR}_+\times M$ with  symplectic form  $\omega = \hd(t\eta)$, metric $\tilde{g}= dt^2+ \eta^{\otimes 2} +tg\mid_V$, and  almost complex structure $\tilde{J}=\hat{J}$.

\medskip

To  answer the first question we use the cylinder case and   a certain natural perturbation of the first canonical connection on $(T\hat{M}, \hg, \hat{J})$.     This new  connection  on $T\hat{M}$ preserves the splitting $T\hat{M}={\bR}\pa_t\oplus TM$ and induces a connection on $TM$ with   the required properties (see \ref{ss: 31}). Moreover, when $M$ is a $CR$ manifold this connection  coincides with the Webster connection, \cite{Stan, Web}.

To answer the second question we use the symplectization $\tilde{M}$ and a   natural perturbation of the Chern connection on $T\tilde{M}$. We obtain a  new connection on $\tilde{M}$ whose restriction to $\{1\}\times M$ is a contact connection (see \ref{ss: 33}). When $M$ is $CR$ this contact connection is also $CR$, but it   never coincides with the Webster connection.  We are not aware whether this contact connection has been studied before.

 These two connections are examples of  geometric connections.  In fact we prove a  much stronger result.

\medskip

\noindent {\bf Theorem.}\hspace{.3cm} (a) {\em On any metric contact  manifold    there exists  a nice  contact connection adapted to $\h$ and a nice contact connection Dirac equivalent to the  Levi-Civita connection. If the manifold is $CR$ these connections are also $CR$.

\noindent (b) On a $CR$ manifold  each Dirac equivalence class of  connections contains at most one nice  $CR$ connection. Moreover, the Webster connection is the unique nice $CR$  connection adapted to $\h$.}

\medskip

Finally, we present several Weitzenb\"{o}ck formul{\ae}  involving   the operator $\h$ (see \ref{ss: 32}). We expect these facts will have applications to  three dimensional Seiberg-Witten theory.

\tableofcontents

\section{General properties  geometric Dirac operators}
\setcounter{equation}{0}

\subsection{Dirac operators  compatible with a metric connection}
\label{ss: 11}

Suppose $(M, g)$ is an oriented, $n$-dimensional  Riemannian manifold.     We will denote a generic local, oriented, synchronous frame  of $TM$ by $(e_i)$.  Its dual coframe  is denoted by $(e^i)$. We will denote the natural duality between a vector space and its dual by $\lan\bullet,\bullet\ran$.

A {\em metric connection} on $TM$ is a connection $\nabla$ on $TM$ such that
\[
X\cdot g(Y,Z)=g(\nabla_XY, Z)+ g(Y, \nabla_XZ),\;\;\forall X,Y,Z\in {\rm Vect}\,(M).
\]
The {\em torsion} of  a metric connection $\nabla$ is the $TM$-valued $2$-form $T=T(\nabla)$  defined by
\[
T(X,Y)=\nabla_XY-\nabla_YX-[X,Y].
\]
The Levi-Civita connection, denoted by $D$ in the sequel is the metric connection uniquely determined by the condition $T(D)=0$. Any  metric connection $\nabla$    can be uniquely written as $D + A$, where $A\in \Omega^1(\End_-(TM))$, where $\End_-$ denotes the space of skew-symmetric endomorphisms. $A$ is called the {\em the potential of $\nabla$}.

There are natural {\em isomorphisms}
\[
\Omega^2(TM)\Lra \Omega^2(T^*M),\;\;T\mapsto T^\dag
\]
\[
\Omega^1(\End_-(TM))\mapsto \Omega^2(T^*M),\;\;A\mapsto A^\dag
\]
defined as follows.
\[
\Omega^2(TM)\ni T \mapsto T^\dag,\;\;\lan X, T^\dag(Y, Z)\ran =g(X, T(Y,Z))
\]
and
\[
\Omega^1(\End_-(TM))\ni A\mapsto A^\dag,\;\; \lan X, A^\dag(Y,Z)\ran = g(A_XY, Z)=:A^\dag(X;Y,Z),
\]
$\forall X,Y,Z\in {\rm Vect}\,(M)$. In local coordinates, if \
\[
T(e_j,e_k)=\sum_{i}T^i_{jk}e_i,\;\;A_{e_i}e_j=\sum_k A^k_{ij}e_k
\]
 then
\[
T^\dag(e_j,e_k)=\sum_iT^i_{jk}e^i,\;\;A^\dag(e_j,e_k)= \sum_i A_{ij}^ke^i,
\]
or equivalently, $T^\dag_{ijk}=T^i_{jk}$, $A^\dag_{ijk}=A^k_{ij}$.  To simplify the exposition, when working in local coordinates,   we will write $A_{ijk}$ instead of $A^\dag_{ijk}$ etc. Define
\[
{\rm tr}:\Omega^2(T^*M)\ra \Omega^1(M),\;\; \Omega^2(T^*M)\ni (B_{ijk})\mapsto ({\rm tr}\,B)=\sum_{i,k}B_{iik}e^k
\]
and  the {\em Bianchi map}
\[
{\gb}:\Omega^2(T^*M)\ra \Omega^3(M),
\]
\[
\Omega^2(T^*M)\ni (B_{ijk})\mapsto {\gb} B=\sum_{i<j<k} (B_{ijk}+B_{kij}+B_{jki})e^i\wedge e^j\wedge e^k.
\]
Note that if $B\in \Omega^3(M)\subset \Omega^2(T^*M)$ then $B=\frac{1}{3}{\gb}B$.

For any  $A\in \End (TM)$ and $\alpha\in \Omega^1(M)$ we define $A\wedge \alpha \in \Omega^2(T^*M)$ by the  equality
\[
(A\wedge \alpha)(X;Y,Z)= (AX)_\flat\wedge \alpha(Y,Z),\;\;\forall X,Y,Z\in \Vect(M),
\]
where ${\bullet}_\flat$ (resp. ${\bullet}^\flat$) denotes the  $g$-dual of a vector (resp. covector) $\bullet$. The following elementary result lists some basic properties of  the above operation.

\begin{lemma}   Let $A\in =\End(TM)$, $\alpha\in \Omega^1(M)$ and  set
\[
A_+=\frac{1}{2}(A+A^*),\;\;A_-=\frac{1}{2}(A-A^*).
\]
Then
\[
\tr (A\wedge \alpha)= (\tr A)\alpha- A^t\alpha,
\]
and
\[
{\gb}(A\wedge \alpha) =2\omega_{A_-}\wedge \alpha,
\]
where
\[
\omega_{A_-}(X,Y)=g(A_-X,Y),\;\;\forall X,Y\in \Vect(M).
\]
\label{lemma: op-wedge}
\end{lemma}
Using the above operations  we can orthogonally decompose $\Omega^2(T^*M)$ as
\[
\Omega^2(T^*M)=\Omega^1(M)\oplus \Omega^3(M)\oplus \Omega^2_0(T^*M)
\]
where
\[
\Omega^2_0:=\Bigl\{ A\in \Omega^2(T^*M);\;\;{\gb}A=\tr A=0\Bigr\},
\]
and $\Omega^1(M)$ embeds in $\Omega^2(T^*M)$  via the map
\[
\Omega^1(M)\ra\Omega^2(T^*M),\;\;\alpha\mapsto\tilde{\alpha}:= \frac{1}{n-1} \Bigl({\bf 1}_{TM}\wedge \alpha\bigr)
\]
Using this orthogonal splitting we can decompose any $A\in \Omega^2(T^*M)$ as
\[
A=\widetilde{\tr A} + \frac{1}{3}{\gb}A + P_0A,\;\; P_0A:=A-\tilde{\tr A}-  \frac{1}{3}{\gb}A \in \Omega^2_0(T^*M).
\]
The next result,  whose  proof can be found in \cite{Gau}, states  that a   metric connection is determined by its torsion in a very explicit way.

\begin{proposition}  Suppose that  $\nabla$ is a metric connection with potential $A$ and torsion $T$. Then
\begin{equation}
T^\dag=-A^\dag+{\gb}A^\dag,
\label{eq: pot-tors}
\end{equation}
\begin{equation}
A^\dag=-T^\dag +\frac{1}{2}{\gb}T^\dag.
\label{eq: tors-pot}
\end{equation}
In particular
\[
{\gb}A^\dag=\frac{1}{2}{\gb}T^\dag,\;\; {\rm tr}\,A^\dag=-{\rm tr}\,T^\dag.
\]
\label{prop: pot-tors}
\end{proposition}
Since all the computations we are about to perform are local we can assume that $M$ is equipped with a $spin$ structure  and we denote by ${\bS}$ the associated complex spinor bundle\footnote{${\bS}$ is ${\bZ}_2$-graded if $n=\dim M$ is even and it is ungraded if $n$ is odd.}.  We have a {\em Clifford multiplication map}
\[
\bc: \Omega^*(M)\ra \End ({\bS}).
\]
A hermitian connection $\tna$ on ${\bS}$ is said to be compatible  with the Clifford multiplication and the  metric connection $\nabla$ on $TM$ if
\[
\tna_X\Bigl(\bc(\alpha) \psi\Bigr)= \bc(\nabla_X\alpha)\psi+\bc(\alpha)\tna_X \psi,\;\;\forall X\in {\rm Vect}\,(M),\;\;\alpha\in \Omega^1(M),\;\;\psi\in C^\infty({\bS}).
\]
We denote by ${\gA}_\nabla={\gA}_\nabla({\bS})$ the  space of hermitian connections on ${\bS}$ compatible with the Clifford multiplication and $\nabla$.
\begin{proposition} The space ${\gA}_\nabla({\bS})$ is an affine space modelled  by the space $\ii \Omega^1(M)$ of imaginary $1$-forms on $M$.
\end{proposition}

\noindent {\bf Proof}\hspace{.3cm}  Suppose $\tna^0, \tna^1\in {\gA}_\nabla$. Set $C:= \tna^1-\tna^0\in \Omega^1({\rm End}\,({\bS})\,)$. Since  both $\tna^i$, $i=0,1$, are compatible with the Clifford multiplication and $\nabla$ we deduce  that for every $X\in {\rm Vect}\,(M)$  the endomorphism $C(X):=X\inpr C$ commutes with the Clifford multiplication. Since the  fibers of ${\bS}$ are irreducible   Clifford modules  we deduce from Schur's Lemma that $C(X)$  is a constant in each  fiber, i.e  $C\in \Omega^1(M)\otimes {\bC}$. Since both $\tna^i$ are hermitian connections we conclude that $C$ must  be purely imaginary $1$-form. $\Box$

\medskip
\begin{definition} A {\em geometric Dirac operator} on ${\bS}$ is a first order  partial differeintial operator $\dir$ of the form
\[
\dir=\dir(\tna):C^\infty({\bS})\stackrel{\tna}{\Lra}C^\infty(T^*M\otimes {\bS})\stackrel{\bc}{\ra} C^\infty({\bS})
\]
where $\tna\in {\gA}_\nabla({\bS})$ for some metric connection $\nabla$ on $TM$. The   geometric Dirac operator is called {\em nice} if it is formally self-adjoint.
\end{definition}
Locally, a geometric Dirac operator has the form
\[
\dir(\tna) =\sum_i\bc(e^i)\tna_{\e_i}.
\]
Every metric connection $\nabla$  canonically determines a connection $\nah\in {\gA}_\nabla({\bS})$ locally described  as follows. If the $\uso(n)$-valued 1-form $\omega$ associated by the frame $(e_i)$ to the connection $\nabla$ is defined by
\[
\nabla e_j =\sum_{i,k}e^k\otimes \omega^i_{kj}e_i,\;\;\omega^i_{kj}+\omega^j_{ki}=0,
\]
then the  induced connection on ${\bS}$ is given by the ${\rm End}_-\,({\bS})$-valued 1-form (see \cite{N1})
\begin{equation}
\hat{\omega}=-\frac{1}{4}\sum_{i,j,k}e^k\otimes \omega^i_{kj}\bc(e^i)\bc(e^j).
\label{eq: nah}
\end{equation}
We set $\dir(\nabla):=\dir(\nah)$ and $\dir_0:=\dir(\hat{D})$. $\dir_0$ is the usual $spin$ Dirac operator. We see that every geometric operator has the form
\[
\dir=\dir(\nabla)+\bc(\ii a)
\]
where $\nabla$ is a metric connection on $M$ and $a\in \Omega^1(M)$.  The connection $\nabla$ is called {\em nice} if $\dir(\nabla)$ is nice. We denote by ${\gA}_{nice}(M)$ the space of nice connections on $M$.

\begin{proposition}  The connection  $\nabla$  with torsion $T$ is nice if and only if $\tr T^\dag=0$.
\label{prop: nice}
\end{proposition}

\noindent {\bf Proof}\hspace{.3cm} Note that
\begin{equation}
\nai e_j =\naj e_i +T_{ij},\;\;\forall i,j
\label{eq: 21}
\end{equation}
and
\begin{equation}
{\bf div}_g(e_i)=0,\;\;\forall i.
\label{eq: 22}
\end{equation}
 We have (at $x_0$)
\[
\dir^* =\sum_k \nah^*_k\bc(e^k)^* =\sum_k \tna_k\bc(e^k)=\sum_k \bc(\nak e^k)+\sum_k \bc(e^k)\nah_k=\bc\left(\sum_k \nak e^k\right) +\dir.
\]
Thus  $\nabla$ is nice if and only  if
\[
\bc\left(\sum_k \nak e^k\right)=0
\]
We compute easily that
\[
(\naj e^i)(e_k)= -e^i(\naj e_k)= -g(e_i,\naj e_k)=-g(e_i, \nak e_j +T_{jk})
\]
so that
\begin{equation}
\naj e^i = -\sum_k g(e_i, \nak e_j +T_{jk}) e^k.
\label{eq: deriv}
\end{equation}
Hence
\[
\sum_k \nak e^k = -\sum_k \sum_i g(e_k, \nai e_k + T_{ki}) e^i
\]
$(g(e_k, \nai e_k) =0$ at $x_0$)
\[
=-\sum_i\left(\sum_k g(e_k,T_{ki}) \right) e^i.
\]
This concludes the proof of the proposition.  $\Box$

\medskip

\begin{proposition} Suppose that $\nabla=D+A$ is a nice connection on $TM$. Then
\[
\dir(\nabla)=\dir_0+ \frac{1}{2}\bc({\gb}A^\dag)=\dir_0+\frac{1}{4}\bc({\gb}T^\dag).
\]
\label{prop: geom-dir1}
\end{proposition}

\noindent {\bf Proof}\hspace{.3cm} Observe that
\[
\nah= \hat{D}-\frac{1}{4}\sum_{i,j,k}e^k\otimes A^i_{kj}\bc(e^i)\bc(e^j)=\hat{D}-\frac{1}{4}\sum_{i,j,k}e^k\otimes A_{kji}\bc(e^i)\bc(e^j)
\]
so that
\[
\dir(\nabla)-\dir_0=-\frac{1}{4}\sum_{i,j,k}A_{kji}\bc(e^k)\bc(e^i)\bc(e^j)=\frac{1}{4}\sum_{i,j,k}A_{kij}\bc(e^k)\bc(e^i)\bc(e^j).
\]
Since $\tr A =0$ we deduce  that    the contributions corresponding to triplets $(i,j,k)$  where two entries are identical add up to zero. Hence
\[
\frac{1}{4}\sum_{i,j,k}A_{kij}\bc(e^k)\bc(e^i)\bc(e^j)=\frac{1}{4}\sum_{i<j<k}\Bigl(\, ({\gb}A)_{ijk}-({\gb}A)_{jik}\Bigr)\bc(e^i)\bc(e^j)\bc(e^k)=\frac{1}{2}\bc({\gb}A).\;\;\Box
\]
\begin{corollary} Suppose $\dir=\dir_0+\bc(\varpi)$, $\varpi\in \Omega^3(M)$. Then
\[
\dir=\dir(\nabla)
\]
where $\nabla= D+ A,\;\; A^\dag= \frac{2}{3}\varpi$. $\Box$
\label{cor: geom-dir1}
\end{corollary}

The above result can also be rephrased in the language of superconnections described  e.g. in \cite{BGV}. Suppose $\varpi\in \Omega^3(M)$. The operator $d  + \bc(\varpi)$ is a superconnection on  the trivial line bundle $\uc$. Taking the tensor product  it with the connection $\hat{D}$ on ${\bS}$ we obtain a superconnection on ${\bS}={\bC}\otimes {\bS}$
\[
{\bA}_\varpi:= \varpi\otimes {\bf 1} +{\bf 1}\otimes \hat{D}:C^\infty({\bS})\ra \Omega^*({\bS}).
\]
The Dirac operator determined by this  superconnection  is
\[
\bc\circ {\bA}_\varpi=\dir_0+\bc(\omega).
\]
Two connections $\nabla^0,\nabla^1\in {\gA}_{nice}(M)$ will be called   {\em Dirac equivalent} if
\[
\dir(\nah^0)=\dir(\nah^1).
\]
The above results  show that two connections $\nabla^0$ and $\nabla^1$ are Dirac equivalent if and only if
\[
\bc({\gb}T(\nabla^0)^\dag)=\bc({\gb}T(\nabla^1)^\dag)\Longleftrightarrow{\gb}T(\nabla^1)^\dag={\gb}T(\nabla^0)^\dag.
\]

\subsection{Weitzenb\"{o}ck formul{\ae}}
\label{ss: 12}

Suppose $(E,h)$ is a Hermitian vector bundle over $M$.  A  {\em generalized Laplacian}
is  a {\em formally self-adjoint}, second order partial
differential operator $L: C^\infty(E)\ra C^\infty(E)$ whose  principal symbol  satisfies
\[
\si_L(\xi)=-|\xi|^2_g {\mathbf 1}_E.
\]
The following classical result   is the basis of all the constructions in this section.   We  include here a  proof because of its relevance in the sequel.

\begin{proposition}{\bf (\cite[Sec. 2.1]{BGV}, \cite[Sec. 4.1.2]{Gky})} Suppose $L$ is a    generalized Laplacian on  $E$. Then there exists a {\em unique hermitian} connection $\tilde{\nabla}$ on $E$ and a {\em unique selfadjoint}  endomorphism ${\cal R}$ of $E$ such that
\begin{equation}
L =\tilde{\nabla}^*\tna +{\cal R}
\label{eq: laplace}
\end{equation}
We will refer to this  presentation of a generalized Laplacian as  {\em  the Weitzenb\"{o}ck presentation} of $L$.
\label{prop: laplace}
\end{proposition}

\noindent{\bf Proof}\hspace{.3cm}  Choose an arbitrary hermitian connection $\nabla$ on $E$.  Then $L_0=\nabla^*\nabla$ is a generalized Laplacian so that $L-L_0$ is a first order operator which can be represented  as
\[
L-L_0 = A \circ \nabla +  B
\]
where
\[
A: C^\infty(T^*M\otimes E) \ra C^\infty(E)
\]
is a  bundle morphism and $B$ is an endomorphism of $E$. We  will regard $A$ as an ${\rm End}\,(E)$-valued 1-form on $M$.  Hence
\begin{equation}
L= \nabla^*\nabla + A\circ \nabla +B.
\label{eq:  lapl1}
\end{equation}

The connection $\nabla$ induces a connection  on ${\rm End}(E)$ which we continue to denote with $\nabla$
\[
\nabla: C^\infty({\rm End}\,(E))\ra \Omega^1({\rm End}\,(E)).
\]
We define the {\em divergence} of $A$ by
\[
{\bf div}_g (A):= -\nabla^* A.
\]
If $(e_i)$ is a local synchronous frame at $x_0$  and, if $A=\sum_i A_ie^i$, then,  at $x_0$, we have
\[
{\bf div}_g(A)=\sum_i\nai A_i.
\]
Note that since $(L-L_0)=\sum_i A_i\nai + B$ is formally selfadjoint we deduce
\begin{equation}
A^*_i=-A_i,  \;\;{\bf div}_g(A)=B-B^*.
\label{eq: 23}
\end{equation}
We seek a hermitian connection $\tna = \nabla + C$ , $C\in \Omega^1({\rm End}\,(E))$  and an endomorphism ${\cal R}$ of $E$ such that
\[
\tna^*\tna +{\cal R} =\nabla^* \nabla + A\circ \nabla  +  B.
\]
We set $C_i:=e_i\inpr C$ so that  we have the local   description
\[
\tna= \sum_i e^i\otimes (\nai + C_i),\;\; C_i^*=-C_i,\;\;\forall i.
\]
Then, as in \cite{N1},  Example 9.1.26, we deduce that, at $x_0$
\[
\tna^*\tna = -\sum_i(\nai+C_i)(\nai+C_i)
\]
$(\lan C_i\ran ^2:= C_iC_i^*=-C_i^2$)
\[
=-\sum_i \nai^2 -\sum_i \nai C_i -2\sum_iC_i\nai + \sum_i \lan C_i\ran ^2
\]
($\lan C\ran^2=\sum_i\lan C_i\ran^2$)
\[
= \nabla^* \nabla  -2C\circ \nabla  - {\bf div}_g(C)+\lan C\ran^2=\nabla^*\nabla +A\circ\nabla +B-{\cal R}.
\]
We deduce immediately that
\begin{equation}
C=-\frac{1}{2} A, \;\; {\cal R}= B -\frac{1}{2}{\bf div}_g (A)
-\lan C\ran^2 \stackrel{(\ref{eq: 23})}{=} \frac{1}{2}(B+B^*)-\frac{1}{4}\lan A\ran ^2.
\label{eq: 24}
\end{equation}
The proposition is proved. $\Box$

\bigskip

If $\dir$ is a geometric Dirac operator on ${\bS}$ then  both $\dir^*\dir$ and $\dir\dir^*$  are generalized Laplacians. Suppose now that $\nabla$ is a nice connection  on our spin manifold $(M,g)$.  It determines a  nice Dirac operator $\dir(\nabla)$. We  denote by $\nabla^{\gw}$ and respectively $\R_\nabla$ the Weitzenb\"{o}ck connection and respectively remainder of the generalized Laplacian $\dir(\nabla)^2$. A classical result of Lichnerowicz states  that if $\nabla$ is the Levi-Civita connection then   $\nabla^{\gw}= \nah$ and ${\cal R}=\frac{s}{4}$, where $s$ is the scalar curvature of the Riemann metric  $g$.  When $\nabla$ is not symmetric the situation is more complicated.
We will present some general formul{\ae} describing  $\nabla^{\gw}$ and $\R$.
\[
\dir^2 =\sum_{i,j}\bc(e^i)\nah_i c(e^j) \nah_j = \sum_{i,j} \bc(e^i)\bc(e^j)\nah_i\nah_j + \sum_{i,j}\bc(e^i)\bc(\nai e^j) \nah_j
\]
\[
= -\sum_i \nah_i^2  +\sum_{i<j}\bc(e^i)\bc(e^j)[\nah_i, \nah_j] + \sum_{i,j}\bc(e^i)\bc(\nai e^j) \naj
\]
\[
= \nah^*\nah  + \sum_{i,j}\bc(e^i)\bc(\nai e^j) \naj + \sum_{i<j}\hat{R}_{ij} \bc(e^i)\bc(e^j)
\]
where
\[
\hat{R}=\sum_{i<j}e^i\wedge e^j\hat{R}_{ij}
\]
denotes the curvature of $\nah$. We need to better understand the quantity $A$ (the coefficient of the first order part of $\dir^2$)  which at $x_0$ is defined as
\[
A= \sum_{i,j}e^j\otimes \bc(e^i)\bc\left(\nai e^j\right).
\]
Using (\ref{eq: deriv}) we deduce
\[
A =\sum_{i,j}e^j\otimes \bc(e^i) \bc(-\sum_k \lan e_j, \nak e_i +T_{ik}\ran e^k)
\]
\[
=-\sum_{i,j,k}e^j\otimes \lan e_j, \nak e_i \ran \bc(e^i)\bc(e^k) -\sum_{i,j,k} e^j\otimes \lan e_j, T_{ik}\ran \bc(e^i)\bc(e^k)
\]

\[
= \sum_je^j\sum_k \lan e_j, \nak e_k\ran -\sum_j e^j \otimes \sum_{i\neq k}\lan e_j, \nak e_i \ran \bc(e^i)\bc(e^k) -\sum_{j, i,k}e^j \otimes \lan e_j ,T_{ik}\ran \bc(e^i)\bc(e^k)
\]
($\lan e_j, \nak e_k\ran=-\lan \nak e_j, e_k\ran$ at $x_0$, $\nabla_k e_i-\nabla_ie_k=T_{ki}=-T_{ik}$)
\[
= -\sum_je^j\sum_k\lan \nak e_j, e_k\ran +\sum_j e^j\otimes \sum_{i<k} \lan e_j, T_{ik}\ran \bc(e^i)\bc(e^k)
\]
\[
 -2 \sum_j e^j\otimes \sum_{i<k} \lan e_j, T_{ik}\ran \bc(e^i)\bc(e^k)
\]
(switch the order of summation in the first term)
\[
=-\sum_k\Bigl(\sum_j\lan \nak e_j, e_k\ran e^j\Bigr)-\sum_j e^j\otimes \sum_{i<k} \lan e_j, T_{ik}\ran \bc(e^i)\bc(e^k)
\]
\[
=\sum_k \nak e^k-\sum_j e^j\otimes \sum_{i<k} \lan e_j, T_{ik}\ran \bc(e^i)\bc(e^k)
\]
($\sum_k\nabla_ke^k=0$)
\[
= -\frac{1}{2}\sum_{i,j,k} e^j\otimes  \lan e_j, T_{ik}\ran \bc(e^i)\bc(e^k) =: -\alpha(T).
\]
We deduce
\begin{equation}
\dir^2 =\nah^* \nah -\alpha(T)\circ \nah + c(\hat{R})
\label{eq: square}
\end{equation}
where
\[
c(\hat{R}):= \sum_{i<j}\bc(e^i)\bc(e^j) \hat{R}_{ij}.
\]
Using the equalities (\ref{eq: 24}) we reach the following conclusion.

\begin{proposition} We have the Weitzenb\"{o}ck formula
\[
\dir^2=(\nabla^{\gw})^*\nabla^{\gw}+\R_\nabla
\]
where
\begin{equation}
{\nabla}^{\gw} = \nah +\frac{1}{2}\alpha(T)=\nah+\frac{1}{4}\sum_{i,j,k}e^i\otimes T_{ijk}\bc(e^j)\bc(e^k).
\label{eq: weiz1}
\end{equation}
\begin{equation}
{\mathcal R}_\nabla= \frac{1}{2}(c(\hat{R}) +c(\hat{R})^*)-\frac{1}{4}\lan \alpha(T)\ran^2,
\label{eq: weiz2}
\end{equation}
where $T$ denotes the  torsion of $\nabla$ and $\hat{R}$ the curvature of $\nah$. $\Box$
\label{prop: weiz}
\end{proposition}

\bigskip

\begin{remark} Observe that $\nabla^{\gw}$ is the connection on ${\bS}$ induced by the nice connection $\nabla' = \nabla +A$ where $A^\dag= T^\dag$. Using (\ref{eq: pot-tors}) we deduce
\[
T^\dag(\nabla')=T(\nabla)^\dag -A^\dag +{\gb}A^\dag={\gb}T(\nabla)^\dag
\]
\end{remark}

The Weitzenb\"{o}ck remainder can be given a more explicit description.   More precisely  we know from Proposition \ref{prop: geom-dir1}  that
\[
\dir(\nabla)=\dir_0+\frac{1}{4}\bc({\gb}T^\dag).
\]
We  set $\varpi:=\frac{1}{4}{\gb}T^\dag$.  As explained at the end of \ref{ss: 11}, $\dir(\nabla)$ is the Dirac operator associated to the superconnection $\hat{D}+\varpi$.   Using \cite[Thm. 1.3]{Bis} we deduce that the Weitzenb\"{o}ck remainder of $\dir^2$ is
\[
\R_\nabla=  \frac{s}{4} + \bc(d\varpi) -2\|\varpi\|^2=\frac{s}{4}+\frac{1}{4}\Bigl(\bc(d{\gb}T^\dag) -\frac{1}{2}\|{\gb}T^\dag\|^2\Bigr).
\]
where $\|\bullet\|$ denotes the  pointwise norm of a differential form and $s$ denotes the scalar curvature of $g$.

The following result summarizes the main facts we proved so far.

\begin{theorem} Denote  by $\dir_{spin}$ the $spin$-Dirac operator induced by the Levi-Civita $D$, $\dir_{spin}=\dir(\hat{D})$. Any geometric Dirac operator $\dir$ can be written as
\[
\dir=\dir_{spin}+ \bc(\varpi) +\bc(\ii a),\;\;a\in \Omega^1(M), \;\;\varpi\in \Omega^3(M).
\]
Additionally,  if $\nabla=D +\frac{2}{3}\varpi+ U$, where $U\in \Omega^2(T^*M)$ is such that
\[
\tr U=0={\gb}U=0
\]
then
\[
\dir= \dir(\nah) + \bc(\ii a)
\]
and
\[
\dir(\nah)^2 = \bigl(\nabla^{\gw}\bigr)^*\nabla^{\gw} +{\mathcal R}_\nabla +\bc(\ii da)
\]
where
\[
{\mathcal R}_\nabla =\frac{1}{4} s(g) +(\bc(d\varpi)-2\|\varpi\|^2)
\]
\label{th: geom-dir}
\end{theorem}
The last theorem   has an obvious extension  where we replace ${\bS}$ by the complex spinor bundle ${\bS}_\si$  determined by a $spin^c$-structure $\si$ on $M$.     This case requires the choice of a hermitian connection on the line bundle $\det {\bS}_\si$. In the spin case $\det {\bS}\cong {\bC}$ and the additional hermitian connection on the trivial line bundle is encoded by the imaginary $1$-from $\ii a$ appearing in the statement of Theorem \ref{th: geom-dir}.

\section{Dirac operators on almost-hermitian manifolds}
\setcounter{equation}{0}

\subsection{Basic differential geometric objects on an  almost-hermitian manifolds}
\label{ss: 21}

In this subsection we  survey a few  differential geometric facts
concerning almost complex manifolds. For more details we refer to
\cite{Gau, KN, Lich} which served  as sources of inspiration.

Consider an almost-hermitian manifold $(M^{2n}, g, J)$. Recall
that this means  that $(M,g)$ is a Riemann manifold and $J$ is a skew-symmetric endomorphism of $TM$ such that $J^2=-1$.  Fix $x_0\in M$  and $(e_1, f_1, \cdots, e_n, f_n)$ a local, oriented orthonormal frame of $TM$. We also assume it is adapted to  $J$ that is
\[
f_j=Je_j,\;\;\forall j=1,\cdots, n.
\]
We denote by $(e^1,f^1,\cdots, e^n, f^n)$ the dual coframe. Let $\ii:=\sqrt{-1}$ and   fix one such adapted local frame.    We split $TM\otimes {\bC}$ into $\pm \ii$-eigen-subbundles of $J$, $TM^{1,0}$ and $T^{0,1}$. These are naturally equipped with hermitian metrics induced by $g$ and have  natural local unitary  frames near $p_0$
\[
TM^{1,0}:\;\;\ve_k:=\frac{1}{\sqrt{2}}(e_k-\ii f_k),\;\;k=1,\cdots, n,
\]
\[
TM^{0,1}:=\;\;\bar{\ve}_k:=\frac{1}{\sqrt{2}}(e_k+\ii f_k),\;\;k=1,\cdots, n.
\]
Form by duality $T^*M^{1,0}$ and $T^*M^{0,1}$ with local unitary frames   given by
\[
\ve^k:=\frac{1}{\sqrt{2}}(e^k+\ii f^k),\;\;k=1,\cdots, n
\]
and respectively,
\[
\bar{\ve}^k:=\frac{1}{\sqrt{2}}(e^k-\ii f^k),\;\;k=1,\cdots, n.
\]
We have unitary decompositions
\[
\Lambda^mT^*M\otimes {\bC}=\bigoplus_{p+q=m} \Lambda^{p,q}T^*M,\;\;m=0,\cdots, 2n
\]
where
\[
\Lambda^{p,q}T^*M:= \Lambda^pT^*M^{1,0}\otimes \Lambda^qT^*M^{0,1}.
\]
Set $K_M:= \Lambda^{n,0}T^*M$. We denote by $P^{p,q}$ the unitary projection onto $\Lambda^{p,q}$ and define
\[
\bpar: \Omega^{p,q}(M)\ra \Omega^{p,q+1}(M),\;\;\bpar:= P^{p,q+1}\circ d
\]
and
\[
\partial:\Omega^{p,q}(M)\ra \Omega^{p+1,q}(M),\;\;\partial:=P^{p+1,q}\circ d.
\]
Define $d^c:\Omega^p(M)\ra \Omega^{p+1}(M)$ by
\[
d^c\alpha(X_0,X_1,\cdots , X_p)=\alpha(-JX_0,-JX_1,\cdots, -JX_p).
\]
The space $\Omega^3(M)\otimes {\bC}$ splits unitarily as
\[
\Omega^3\otimes {\bC}= \Omega^+\oplus \Omega^-,
\]
where
\[
\Omega^+:=\Omega^{2,1}\oplus \Omega^{1,2},\;\;\Omega^-:=\Omega^{3,0}\oplus \Omega^{0,3}.
\]
Finally, introduce the involution ${\gM}$ on $\Omega^2(T^*M)$ defined by
\[
{\gM}B(X;Y,Z)=B(X;JY,JZ).
\]
Observe that
\[
\psi^+={\gb}{\gM}\psi^+,\;\;\forall \psi^+\in \Omega^+.
\]
We denote by $\Omega^{1,1}(T^*M)$ the $1$-eigenspace of ${\gM}$  and by $\Omega^{1,1}_s(T^*M)$ the  intersection of $\ker{\gb}$ to $\Omega^{1,1}(T^*M)$. Thus
\[
A\in \Omega_s^{1,1}(T^*M)\Longleftrightarrow A= {\gM}A,\;\;{\gb}A=0.
\]

The Nijenhuis  tensor $N\in \Omega^2(TM)$ is  defined by
\[
N(X,Y):=\frac{1}{4}([JX,JY]-[X,Y]-J[X,JY]-J[JX,Y]),\;\;\forall X,Y\in {\rm Vect}\,(M).
\]
Notice that $N(JX,Y)= N(X,JY)=-JN(X,Y)$. This implies immediately that $\tr N^\dag =0$.

We denote by $D$ the Levi-Civita connection determined by the metric $g$ and  by $\omega$ the fundamental two form
defined by
\[
\omega(X,Y)=g(JX,Y),\;\;\forall X,Y\in {\rm Vect}\,(M).
\]
Locally we have
\[
\omega=\ii\sum_j\ve^j\wedge \bar{\ve}^j.
\]
The {\em Lee form} $\theta$ determined by $(g,J)$ is defined by
\[
\theta= \Lambda (d\omega)=-J \Lambda\bigl((d^c\omega)^+\bigr),
\]
where $\Lambda$ denotes  the contraction by $\omega$, $\Lambda =(\omega\wedge \;)^*$, and $J$  acts on the $1$-form $\alpha$  by
\[
J\alpha(X)= -\alpha(JX),\;\;\forall X\in {\rm Vect}\,(M).
\]
We have the following identity
\begin{equation}
g((D_XJ)Y,Z)=-\frac{1}{2}d\omega(X,JY,JZ)+\frac{1}{2}d\omega(X,Y,Z) +2g(N(Y,Z),JX).
\label{eq: kn}
\end{equation}
The form $\omega$ determines the skew-symmetric part of $N^\dag$ via the identity
\[
{\gb}N^\dag=(d^c\omega)^-.
\]
The almost complex structure defines a Cauchy-Riemann operator
\[
\bpar_J: C^\infty (TM^{1,0})\ra \Omega^{0,1}(TM^{1,0})
\]
defined by
\[
X\inpr \bpar_J Y= [X,Y]^{1,0},\;\;\forall X\in C^\infty(TM^{0,1}),\;\;Y\in C^\infty(TM^{1,0}).
\]

A {\em Hermitian connection} on $TM$  is a metric connection $\nabla$ such that $\nabla J=0$. A Hermitian connection $\nabla$ is completely determined $\psi_+:=\frac{1}{3}({\gb}T^\dag)^+$ and $B:=(T^\dag)_s^{1,1}$ via the equality (see \cite[Sec. 2.3]{Gau})
\[
T(\nabla)^\dag= N^\dag+\frac{1}{8}(d^c\omega)^+ -\frac{3}{8}{\gM}(d^c\omega^+)+ \frac{9}{8}\psi^+-\frac{3}{8}{\gM}\psi^+ +B.
\]
We will denote the above connection  by $\nabla(\psi^+, B)$. When $B=0$  we write $\nabla(\psi^+)$ instead of $\nabla(\psi^+, B)$. Observe that if $T$ is the torsion of $\nabla(\psi^+, B)$ then
\[
{\gb}T^\dag= {\gb} N^\dag +3\psi^+= (d^c\omega)^- + 3\psi^+= {\gb}N^\dag +3\psi^+.
\]
Using the formul{\ae} \cite[(1.3.5), (1.4.9)]{Gau} and the equality $\psi^+={\gb}{\gM}\psi^+$, $\forall \psi^+\in \Omega^+$   we deduce  that
\[
\tr {\gM}\psi^+= -2J\Lambda \psi^+,\;\;\forall \psi^+\in \Omega^+(M).
\]
Since $\tr N^\dag =0$ we deduce that the trace of the torsion of $\nabla(\psi^+, B)$
\[
\tr T\bigl(\nabla(\psi^+,B)\,\bigr)= \tr B + \frac{3}{4}J\Lambda\Bigl( \,(d^c\omega)^++\psi^+\Bigr)= \tr B-\frac{3}{4}\theta +\frac{3}{4}J\Lambda\psi^+.
\]
\begin{ex} The  {\em first canonical connection} (see \cite[Sec. 2.5]{Gau} or \cite{Lich}) is the  Hermitian connection $\nabla^0$
defined by  $B=0$ and
\[
{\gb}T_0^\dag=(d^c\omega)^- -(d^c\omega^+)
\]
so that $\psi^+ =-\frac{1}{3}(d^c\omega)^+$. Its torsion is
\[
T_0^\dag= N^\dag-\frac{1}{4}\bigl( (d^c\omega)^+ + {\gM}(d^c\omega)^+\bigr).
\]
In general,  it is not a nice connection since $\tr T_0^\dag =-\frac{1}{2}\theta$.
\label{ex: 1cancon}
\end{ex}

\begin{ex} The {\em Chern connection} or the {\em second fundamental connection}, \cite{Gau, Lich}, is the {\em unique}  Hermitian connection $\nabla$ on $TM$ such that
\[
\nabla^{0,1}=\bpar_J.
\]
We will denote it by $\nabla^c$. Alternatively (see \cite[Sec. 2.5]{Gau}), it is the
hermitian connection  defined by $B=0$ and ${\gb}T^\dag=(d^c\omega)^- +(d^c\omega)^+$, i.e it is
determined by  $\psi^+ =\frac{1}{3}(d^c\omega)^+$. Its torsion is given by
\[
T_c^\dag= N^\dag + \frac{1}{2}\Bigl((d^c\omega)^+ -{\gM}(d^c\omega)^+\Bigr).
\]
In general, it is not a nice connection since $\tr T_c^\dag =-\theta$.
\label{ex: 2cancon}
\end{ex}

\subsection{The Hodge-Dolbeault operator}
\label{ss: 22}

The  almost hermitian manifold $M$ is equipped with a canonical $spin^c$ structure and  the associated complex spinor bundle is
\[
{\bS}_c:=\Lambda^{0,*}T^*M=\bigoplus_{p\geq 0}\Lambda^{0,p}T^*M.
\]
Note that $\det {\bS}_c=K_M^{-1}$.  The Chern connection  induces a  hermitian connection $\det \nabla^c$ on $K_M^{-1}$ and we denote  by  $\dir_c$ the  geometric Dirac operator induced by the Levi-Civita connection $D$ and the connection $\det\nabla^c$.

If $M$ is spinnable, then a choice of spin  structure is equivalent to a choice of a square root of $K_M$ and in this case ${\bS}_c:={\bS}\otimes K_M^{-1/2}$.

The bundle ${\bS}_c$  has a natural Dirac type operator, the Hodge-Dolbeault operator
\[
{\h}_J:= \sqrt{2}(\bpar +\bpar^*): C^\infty({\bS}_J)\ra C^\infty({\bS}_J).
\]
We have the following result \cite[Thm.2.2]{Bis} and \cite[Sec.3.6]{Gau}.
\[
{\h}_J = \dir_c -\frac{1}{4}\Bigl\{\bc\bigl((d^c\omega)^+\bigr)-\bc\bigl((d^c\omega)^-\bigr)\Bigr\}.
\]
Using Theorem \ref{th: geom-dir} we deduce that ${\h}_J$ is a geometric Dirac operator, more precisely  ${\gH}_J$ is induced by  $\widehat{\nabla}\otimes {\bf 1}+ {\bf 1}\otimes \det \nabla^c$, where $\nabla$ is the connection
\[
\nabla= D-\frac{1}{6}((d^c\omega)^+-(d^c\omega)^-)
\]
with torsion
\[
T^\dag=\frac{1}{3}(d^c(\omega)^- - (d^c\omega)^+).
\]
A stronger result is true. Using the results in the previous subsection we deduce the following result.

\begin{theorem} For every $B\in \Omega_s^{1,1}(T^*M)$  such that  $\tr B=\frac{1}{2}\theta$ there exists  a   {\bf Hermitian} connection $\nabla^b=\nabla^b(B)$  uniquely determined by the following conditions.

\medskip

\noindent (i) $\nabla^b$ is nice.

\noindent (ii) $\nabla^b$ is  Dirac equivalent to $\nabla^0$.

\label{th: Bis}
\end{theorem}

\noindent {\bf Proof} \hspace{.3cm} Since $\nabla^b=\nabla(\psi^+, B)$  is strongly Dirac equivalent to $\nabla$ we deduce that its torsion satisfies
\[
{\gb}T_b= (d^c\omega)^- -(d^c\omega)^+.
\]
Thus we need to choose $\psi^+=-\frac{1}{3}(d^c\omega)^+$.   Now observe that
\[
0=\tr T_b^\dag = \tr B -\frac{1}{2}\theta=0,\;\;\Box
\]

\begin{definition} We will refer to  any of the connections $\nabla^b$  constructed in Theorem \ref{th: Bis} as a {\bf basic connection} determined by an almost Hermitian structure.
\end{definition}
The torsion of  a basic connection $\nabla^b(B)$ is
\[
T_b^\dag= N^\dag -\frac{1}{4}\Bigl( (d^c\omega)^++ {\gM}(d^c\omega)^+\Bigr) +B.
\]
Observe also  that     the first and second fundamental connection  coincide of  an almost K\"{a}hler  structure coincide and they are both basic.   They  are  precisely  the  connections   used by Taubes, \cite{Taubes}, to analyze   the Seiberg-Witten monopoles on  a symplectic manifold.

For any basic connection $\nabla^b$ we have the following identities (\cite[Sec. 3.5]{Gau})
\begin{subequations}
\begin{equation}
(\bpar\phi)(Z_0,Z_1,\cdots ,Z_p)=\sum_{j=0}^p(-1)^j\nabla^b_{Z_j}\phi(Z_0,\cdots, \hat{Z}_j,\cdots, Z_p),
\label{eq: bparn}
\end{equation}
\begin{equation}
\bpar^*\phi(Z_1,\cdots, Z_{p-1})=-\sum_{i=1}^n\Bigl(e_i\inpr\nabla^b_{e_i}\phi  +f_j\inpr \nabla^b_{f_i}\phi \Bigr)(Z_1,\cdots, Z_{p-1}),
\label{eq: bparna}
\end{equation}
\end{subequations}
$\forall Z_0,\cdots , Z_p\in C^\infty(T^{0,1}M),\;\;\phi\in \Omega^{0,p}(M)$.

\section{Dirac operators on contact  $3$-manifolds}
\setcounter{equation}{0}

\subsection{Differential objects on metric contact manifolds}
\label{ss: 31}

We review  a few basic geometric facts concerning   metric  contact manifolds. For more details we refer to \cite{Blair, Tanno}.

A {\em metric contact manifold} (m.c. manifold for brevity) is an oriented manifold of odd dimension $2n+1$ equipped with a Riemann metric $g$ and a $1$-form $\eta$ such that

\bigskip

\noindent $\bullet$ $|\eta(x)|_g=1$, $\forall x\in M$. Denote by $\xi\in {\rm Vect}\,(M)$ the metric dual of $\eta$ and set $V:= \ker \eta \subset TM$.    $V$ is a hyperplane sub-bundle of $TM$ and we denote by $P_V$ the orthogonal projection onto $V$.

\medskip

\noindent $\bullet$ There exists  $J:TM\ra TM$ such that
\[
d\eta(X,Y)= g(JX, Y),\;\;\forall X,Y\in {\rm Vect}\,(M).
\]
and
\[
J^2X= -X+\eta(X)\xi,\;\;\forall X\in {\rm Vect}\,(M).
\]
\begin{definition} A {\bf   contact metric connection} on $(M^{2n+1}, \eta, J, g)$ is a metric conection such that $\nabla J=0= \nabla\xi$.
\end{definition}

The manifold   $M$ is called  {\em positively oriented} if the orientation induced by  the nowhere vanishing $(2n+1)$-form  $\eta\wedge (d\eta)^n$ coincides with the given orientation of $M$. In this case
\[
dv_g=\frac{1}{n!}\eta\wedge (d\eta)^n
\]
Set $\omega:=d\eta$. The metric $g$ is completely determined by $\eta$ and $J$ via the  equality
\[
g(X,Y) =\eta(X)\eta(Y) + d\eta(X, JY)=\eta(X)\eta(Y) + \omega(X, JY).
\]
We have decompositions
\[
V\otimes {\bC}=V^{1,0}\oplus V^{0,1},\;\;V^*\otimes {\bC}= (V^*)^{1,0}\oplus (V^*)^{0,1}
\]
and we set
\[
K_M:= \det (V^*)^{1,0}.
\]
Set $\Phi:= L_\xi J$.  The operator $\Phi$ is a traceless, symmetric endomorphism of $V$ (see \cite{Blair}).  Since $L_\xi(J^2)=0$ we deduce
\begin{equation}
J\Phi+\Phi J=0  \Longrightarrow (J\Phi)^*=(J\Phi)
\label{eq: phi}
\end{equation}
Define the Nijenhuis tensor  $N\in \Omega^2(TM)$ by
\[
N(X,Y)= \frac{1}{2}\Big\{ J^2[X,Y]+[JX,JY]-J[X,JY]-J[JX,Y]\Bigr\}.
\]
Notice that
\[
N(\xi,X)=-\frac{1}{2}J\Phi X,\;\;\forall X\in \Vect(M).
\]
$(M,g,\eta)$ is a {\em Cauchy-Riemann  manifold} ($CR$ for brevity)   if and only if $JN(X,Y)= 0$, $\forall X,Y\in C^\infty(V)$. Equivalently, this means,
and
\[
N(X,Y)+\omega(X,Y)\xi=-J^2N(X,Y)=0,\;\;\forall X,Y\in C^\infty(V).
\]
In this case, the Nijenhuis tensor can be given the more compact  description
\[
N^\dag = \frac{1}{2} J\Phi \wedge \eta -\eta \otimes d\eta.
\]
In particular, $M$ is a $CR$ manifold when $\dim M=3$. Arguing exactly  as in \cite[p.53]{Blair}  we obtain the following result.

\begin{proposition} If $D$ denotes the Levi-Civita connection  of $(M,g)$ then
\[
g\bigl(\,(D_XJ)Y,Z\bigr)=g\bigl(\, JX, N(Y,Z)\,\bigr)+\frac{1}{2}(\eta\wedge d\eta)(JX,Y,Z).
\]
$\forall X,Y,Z\in {\rm Vect}\,(M)$.
\label{prop: DJ}
\end{proposition}

\bigskip

To each metric contact manifold $M$ we can associate an almost Hermitian manifold $(\hat{M}, \hat{g}, \hat{J})$ defined as follows.
\[
\hat{M}= {\bR}\times M,\;\;\hg=dt^2+ g,\;\;\hat{J}\partial_t=\xi,
\]
We will denote by $\hd$ the exterior differentiation on $\hat{M}$. If we set
\[
\hat{\omega}(X,Y)=\hg(\hat{J}X,Y),\;\;\forall X,Y\in {\rm Vect}\,(\hat{M})
\]
then $\hat{\omega}=dt\wedge \eta +\omega$ and $\hat{d}\hat{\omega}=-dt\wedge \omega$. We deduce that the Lee form $\theta=\Lambda(-dt\wedge d\eta)$  is $-ndt$. We will work with local, oriented orthonormal frames $(e_0,f_0, e_1,\cdots, e_n,f_n)$  adapted to $\hat{J}$ such that
\[
e_0=\pa_t,\;\;f_0=\xi,\;\;e^0=dt\;\;f^0=\eta
\]
\[
\hat{\omega}=\ii\ve^0\wedge \bar{\ve}^0+\ii \sum_{k=1}^n\ve^k\wedge \bar{\ve}^k,\;\;\hat{d}\hat{\omega}= -\frac{\ii}{\sqrt{2}}(\ve^0+\bar{\ve}^0)\wedge \sum_{k=1}^n\ve^k\wedge \bar{\ve}^k.
\]
Hence
\[
\hd^c\hat{\omega}=-\frac{1}{\sqrt{2}}(\ve^0-\bar{\ve}^0)\wedge \sum_{k=1}^n\ve^k\wedge \bar{\ve}^k=-\eta\wedge d\eta
\]
so that $({\gb}\hN^\dag)=(\hd^c\hat{\omega})^-=0$. We have the following result, \cite{Blair}.

\begin{proposition}
\[
\hat{N}(X,Y)= \frac{1}{2}N(X,Y)+\frac{1}{2}\omega(X,Y)\xi,\;\;\forall X,Y\in {\rm Vect}\,(M),
\]
\[
\hat{N}(\pa_t,X)=\frac{1}{4}\Phi X,\;\;\forall X\in {\rm Vect}\,(M).
\]
\end{proposition}

Observe that $\hN^\dag\mid_M= \frac{1}{2}N^\dag+\frac{1}{2}\eta\otimes d\eta$ so that
\[
0={\gb}\hN^\dag\mid_M =\frac{1}{2}{\gb}N^\dag+\frac{1}{2}{\gb}(\eta\otimes d\eta)=\frac{1}{2}{\gb}N^\dag+\frac{1}{2}\eta\wedge d\eta.
\]
Hence
\[
{\gb}N^\dag=-\eta\wedge d\eta.
\]
We want to find $B\in \Omega^{1,1}_s(T^*\hat{M})$ such that $\tr B=-\frac{n}{2}dt$ and the basic connection  it induces on  $T^*\hat{M}$  is compatible with the splitting $\pa_t\oplus TM$. The torsion of such a connection is
\[
\hat{T}_b^\dag= \hat{N}^\dag -\frac{1}{4}\Bigl( (\hd^c\hat{\omega})^+ + {\gM}(\hd^c\hat{\omega})^+\Bigr)+B
\]
\[
= \hat{N}^\dag +\frac{1}{4}( \eta\wedge \omega + {\gM}(\eta\wedge \omega)\Bigr)+B.
\]
Thus ${\gb}T_b^\dag= \eta\wedge d\eta$. Using Proposition \ref{prop: pot-tors} we deduce  that $\nabla^b= D+A$ where
\[
A_b^\dag=\frac{1}{2}{\gb}T_b^\dag -T_b^\dag =\frac{1}{4}\Bigl(\eta\wedge d\eta -{\gM}(\eta\wedge d\eta)\Bigr)-\hat{N}^\dag-B.
\]
Thus, for all $X,Y\in {\rm Vect}\,(M)$ which are $t$-independent we have
\[
\hg(\nabla^b_tX,Y)= A_b^\dag(\pa_t; X,Y)
\]
Since
\[
B(\pa_t;\bullet,\bullet)=0\;\;{\rm and}\;\; \hg(\hat{N}(X,Y),\pa_t)=0,\;\;\forall X,Y\in {\rm Vect}\,(M).
\]
we deduce
\[
\hg(\nabla^b_tX,Y)=-\frac{1}{4}{\gM}(\eta\wedge d\eta)(\pa_t; X,Y)=0.
\]
Similarly, we deduce
\[
\hg(\nabla^b_tX,\pa_t)=A_b^\dag(\pa_t; X,\pa_t)=0.
\]
Thus
\[
\nabla^b_tZ=0,\;\;\forall Z\in {\rm Vect}\,(M).
\]
Since $\nabla^b$ is a metric  connection  we deduce
\[
\hg(\nabla^b_\bullet\pa_t,\pa_t)=0.
\]
On the other hand, $\forall X,Y\in \Vect (M)$ we have
\[
\hg(\nabla^b_X\pa_t,Y)=A_b^\dag(X;\pa_t,Y)
\]
\[
=-\frac{1}{4}{\gM}\eta\wedge d\eta(X,\pa_t,Y)-\hg(\hN(\pa_t,Y),X) -B(X;\pa_t,Y)
\]
\[
=\frac{1}{4}g(X_V,Y_V)-\frac{1}{4}g(\Phi Y, X) -B(X;\pa_t, Y),
\]
where $X_V=P_VX$, $Y=P_VY$. Next, $\forall X,Y\in \Vect(M)$, we have
\[
\hg(\nabla^b_XY,\pa_t)=A^\dag_b(X;Y,\pa_t)=-\frac{1}{4}{\gM}\eta \wedge d\eta(X;Y,\pa_t)-\hg(\hat{N}(Y,\pa_t), X)-B(X;Y,\pa_t)
\]
\[
=-\frac{1}{4}g(X_V,Y_V)+\frac{1}{4}g\bigl(\Phi Y,X\bigr)- B(X;Y,\pa_t).
\]
\begin{lemma} There exists $B_0\in \Omega^{1,1}_s(T^*\hat{M})$ such that $\tr B=-\frac{n}{2}dt$ and
\begin{subequations}
\begin{equation}
B(\pa_t;\bullet,\bullet)=0.
\label{eq: b0a}
\end{equation}
\begin{equation}
B(X;Y,\pa_t)=\frac{1}{4}g(X,\Phi Y)-\frac{1}{4}g(X_V,Y_V),\;\;\forall X,Y\in {\rm Vect}\,(V).
\label{eq: b0b}
\end{equation}
\end{subequations}
\label{lemma: b0}
\end{lemma}

\noindent{\bf Proof}\hspace{.3cm}  Define
\[
B=\frac{1}{4} (\Phi \wedge dt +J\Phi\wedge \eta) -\frac{1}{4}(P_V\wedge dt + JP_V\wedge \eta)+\frac{1}{2}\eta\otimes d\eta
\]
and we set
\[
B_0=\frac{1}{4} (\Phi \wedge dt +J\Phi\wedge \eta),\;\;B_1=-\frac{1}{4}(P_V\wedge dt + JP_V\wedge \eta).
\]
We need to show that this definition is correct,  i.e. the above $B$ satisfies all the  required conditions (\ref{eq: b0a}), (\ref{eq: b0b}) and
\[
\tr B=-\frac{n}{2}dt,\;\;{\gb}B=0
\]
\[
B\in \Omega^{1,1}(T^*M).
\]
Here the elementary properties in Lemma \ref{lemma: op-wedge} will come in handy. Since $\Phi$ and $J\Phi$ are symmetric and traceless we deduce that
\[
\tr B_0=0 ,\;\;{\gb}B_0=0.
\]
The condition $B_0\in \Omega^{1,1}$ follows from the identity $\phi J=-J\Phi$. Now observe that $B_1\in \Omega^{1,1}$ and
\[
{\gb}B_1=-\frac{1}{2}\eta\wedge d\eta,\;\; \tr B_1=-\frac{n}{2}dt.
\]
Finally $\eta\otimes d\eta \in \Omega^{1,1}$, it is traceless and
\[
{\gb}(\eta\otimes d\eta)=\eta \wedge d\eta.
\]
The condition (\ref{eq: b0b}) follows by direct computation. The Lemma follows putting together the above facts. $\Box$

\bigskip

If we choose $B$ as in Lemma \ref{lemma: b0} we deduce
\[
\hg(\nabla^b_\bullet X,\pa_t)=,\;\;\forall X\in {\rm Vect}\,(M).
\]

The above computations show that the basic connection $\nabla^b$  of $(\hat{M},\hg, \hat{J})$ determined by $B_0$  preserves the orthogonal splitting $T\hat{M}=\lan \pa_t\ran\oplus TM$ and thus induces a {\em nice}  contact metric connection $\nabla^w$ on $TM$.  We will  call $\nabla^w$ the {\em generalized Webster connection} of  $M$ for reasons  which will be explained   below.  To compute its torsion observe that
\[
\hat{N}^\dag\mid_M =\frac{1}{2}\Bigl\{ N^\dag + \eta\otimes d\eta\Bigr\},
\]
and ${\gM}(\eta\wedge d\eta)\mid_M=\eta\otimes d\eta$.  Finally
\[
B\mid_M= \frac{1}{4}(J\Phi)\wedge \eta-\frac{1}{4}JP_V\wedge \eta +\frac{1}{2}\eta\otimes d\eta.
\]
Since  on $M$ we have the equality $JP_V=J$, the torsion $T_w$  of $\nabla^w$  given by
\begin{equation}
T_w^\dag= \frac{1}{2}N^\dag+\frac{5}{4}\eta\otimes d\eta +\frac{1}{4}\eta\wedge d\eta +\frac{1}{4}(J\Phi -J)\wedge \eta
\label{eq: tw}
\end{equation}
Moreover,  ${\gb}T_w=\eta\wedge d\eta$.

Suppose now that $M$ is a $CR$-manifold. Then
\[
N^\dag = \frac{1}{2}J\Phi\wedge \eta -\eta\otimes d \eta
\]
and thus
\[
T_w^\dag =\frac{3}{4}\eta\otimes d\eta +\frac{1}{4}\eta\wedge d\eta -\frac{1}{4}(J\wedge \eta) +\frac{1}{2}J\Phi\wedge \eta.
\]
We deduce
\[
T_w(X,Y)=d\eta(X,Y)\xi,\;\;\forall X,Y\in {\rm Vect}\,(V).
\]
In particular, because the distribution $V^{1,0}$ is integrable we deduce
\[
T_w(X,Y)=0,\;\;\forall X,Y\in C^\infty(V^{1,0}).
\]
A contact metric connection with the above property will be called a $CR$ metric connection.  Next observe that for $X,Y\in C^\infty(V)$ we have
\[
g(X,T_w(\xi,Y))=T^\dag(X,\xi,Y)=-\frac{1}{4}d\eta(X,Y) +\frac{1}{4}g(JX,Y) +\frac{1}{2}g(J\Phi X,Y).
\]
Hence
\[
T_w(\xi,Y)=\frac{1}{2} J\Phi Y.
\]
Since $\Phi J= -J\Phi$ we deduce
\[
 JT_w(\xi, X) = - T_w(\xi, JX)
\]
Using \cite[Prop. 3.1]{Tanno},  we deduce that when $M$ is a Cauchy-Riemann manifold, the connection $\nabla^w$ on $(V, J)$ is the Tannaka-Webster connection   determined by the $CR$ structure (see \cite{CH, Stan, Tanno, Web} for more details).  The generalized  Webster connection we have constructed does not agree with the generalized Tannaka connection constructed by S.Tanno in \cite{Tanno} because that  connection  is  not compatible with $J$ if $M$ is not a $CR$-manifold.

Finally let us point out that when $M$ is a $CR$ manifold then
\[
g(\nabla^w_\xi X, Y) = g(D_\xi X, Y)+\frac{1}{2}{\gb}T_w^\dag(\xi,X,Y)-T^\dag_w(\xi; X, Y)=g(D_\xi X -\frac{1}{2}JX,Y)
\]
so that
\[
\nabla^w_\xi= D^V_\xi:= P_VD_\xi-\frac{1}{2}J.
\]

\begin{ex} We consider in great detail the special  case of a metric, contact,  $spin$ $3$-manifold $M$. $M$ is automatically a $CR$-manifold so that the torsion of the  (geberalized) Webster connection satisfies
\[
T_w(X,Y)=\frac{1}{2}d\eta(X,Y)\xi,\;\;T_w(\xi,X)=\frac{1}{2}J\Phi X,\;\;\forall X,Y\in C^\infty(V)
\]
\[
{\gb}T_w^\dag =\eta\wedge d\eta.
\]
The  $spin$ Dirac operator $\dir_0$  on $M$ is related to the Dirac operator $\dir(\nabla^w)$ by the equality
\[
\dir(\nabla^w)=\dir_0 +\frac{1}{4}\bc(\gb T^\dag)=\dir_0 +\frac{1}{4}\bc(\eta\wedge d\eta)=\dir_0-\frac{1}{4}.
\]
When $M$ is Sasakian, i.e. $\Phi =0$,  the above equality shows that $\dir(\nabla^w)$ coincides with the  adiabatic Dirac operator introduced in \cite{N2} (see in particular \cite[Eq.(2.20)]{N2} with $\lambda=\frac{1}{2}$, $\delta=1$).
\end{ex}

Later on   we will need to  compare  the connections $\det \nabla^c$ and $\det \nabla^b$ induced  by the Chern connection $\nabla^c$  and respectively $\nabla^b$ on $K^{-1}_{\hat{M}}$.

\begin{proposition}
\[
\det \nabla^c =\det \nabla^b +\frac{n\ii}{2}\eta.
\]
\label{prop: compare}
\end{proposition}

\noindent {\bf Proof}\hspace{.3cm}  Denote by $\nabla^0$ the first fundamental connection of $(\hat{M}, \hat{J})$. We have
\[
\nabla^b = \nabla^0  -B,
\]
where $B$ is described in Lemma \ref{lemma: b0}. Set $\delta:= \ve_0\wedge \ve_1\wedge \cdots \wedge \ve_n$. Then for every vector field  $X$ on $\hat{M}$ we have
\[
\det \nabla_X^b \delta  =\det \nabla_X^0 \delta  - B_X\delta
\]
Observe that
\[
B_X\ve^k = \sum_{j=0}^n C_k^j \ve_j
\]
so that $B_X\delta =(\sum_{j=0}^n C_k^k)\delta$.   On the other hand, $C_k^k= g_c(B_X\ve_k,\bar{\ve}_k)$ where $g_c$ denotes the       complex bilinear extension of $g$.
\[
C_k^k=\frac{1}{2}g_c\bigl(B_X (e_k-\ii f_k)\, , \, e_k +\ii f_k\bigr)=\ii g(B_Xe_k,Je_k)+\ii g(B_Xf_k, Jf_k)
\]
Thus
\begin{equation}
\sum_k C_k^k = -\ii \sum_{k=0}^n \Bigl(\, g(JB_Xe_k,e_k)+g(JB_Xf_k, f_k)\,\Bigr)=-\ii\tr JB_X.
\label{eq: sum}
\end{equation}
The equality
\[
B= \frac{1}{4} \Bigl\{ (\Phi +P_V) \wedge dt -(J\Phi +JP_V)\wedge \eta\Bigr\} +\frac{1}{2}\eta\otimes d\eta
\]
so that
\[
\hg(B_XY, JY)=\frac{1}{4}\Bigl\{\hg(\Phi X,Y)dt(JY) -\hg(J\Phi X,Y)\eta(JY)\Bigr\}
\]
\[
+\frac{1}{4}\Bigl\{ \hg(P_VX,Y)dt(JY) -\hg(JP_VX,Y)\eta(JY)\Bigr\} +\frac{1}{2}\eta(X)d\eta(Y,JY).
\]
We see that $\tr JB_X \neq 0$ only if $X=\xi$ in which case shows that the sum (\ref{eq: sum}) is  $n$.  Hence
\[
\nabla^b\delta =\nabla^0\delta -\ii n \eta.
\]
On the other hand we have the identity, \cite[Eq. (2.7.6)]{Gau},
\[
\det \nabla^c  =\det \nabla^0 +\frac{\ii}{2} J\theta =\det \nabla^0 -\frac{n\ii}{2} J dt= \det \nabla^b+\frac{n\ii}{2}\eta. \;\;\Box
\]
\begin{corollary}
\[
F(\det\nabla^c)=F(\det\nabla^b)+\frac{n\ii}{2} d\eta.\;\;\;\Box
\]
\label{cor: compare}
\end{corollary}

\subsection{Geometric Dirac operators on contact manifolds}
\label{ss: 32}

Consider the Hodge-Dolbeault operator  $\hat{\h}$ on $\hat{M}$
\[
\hat{\h}=\sqrt{2}(\bpar+\bpar^*): \Omega^{0,*}(\hat{M})\ra \Omega^{0,*}(\hat{M}).
\]
It is a geometric Dirac operator  and it is
\[
\hat{\h}=\sqrt{2}\sum_{k=0}^n (\hbc(\ve^k)\nah_{\ve_k}+\hbc(\bar{\ve})\nah_{\bar{\ve}_k})
\]
where  $\hbc$ denotes the Clifford multiplication on $\hat{\bS}_c\cong \Lambda^{0,*}T^*\hat{M}$, $\nah=\nah^b\otimes {\bf 1}+{\bf 1}\otimes \det \nabla^c$, and  $\det \nabla^c$ denotes the Hermitian connection on $K_{\hat{M}}^{-1}$ induced by the Chern connection on $\hat{TM}$. More precisely
\[
\hbc(\bar{\ve}^k)=\sqrt{2}\bar{\ve}^k\wedge\bullet,\;\;\hbc(\ve^k)= -\sqrt{2}\ve^k\inpr \bullet\;.
\]
Above, $\ve^k\inpr\bullet$ denotes the odd derivation of $\Omega^{0,*}(\hat{M})$ uniquely determined by the requirements
\[
\ve^k\inpr\bar{\ve}^j=\delta_{kj},\;\;\forall j,k=0,\cdots n.
\]
We want to point out that
\[
(\bar{\ve}^k\wedge)^*= \ve^k\inpr.
\]
 We set
\[
\J:=\hbc(dt)= \frac{1}{\sqrt{2}}\hbc(\ve^0)+\hbc(\bar{\ve}^0)),\;\;{\bS}_c:=\hat{\bS}^+_c\mid_{0\times M}.
\]
Note that
\[
\hat{\bS}_c\mid_M\cong {\bS}_c\oplus \J {\bS}_c.
\]
The  metric contact structure on $M$  produces a $U(n)$-reduction  of the tangent bundle $TM$ which  in general has only a $SO(2n+1)$-structure. This $U(n)$-reduction induces a $spin^c$ structure on $M$ and ${\bS}_c$ is the associated bundle of complex spinors and
\[
\det{\bS}_c\cong K_M^{-1}.
\]
 The Clifford multiplication on ${\bS}_c$ is  defined by the equality
\[
\bc(\alpha)= \J\hbc(\alpha),\;\;\forall \alpha\in \Omega^1(M).
\]
Along $M$ we can identify $\hat{\bS}_c^-$ with $\J{\bS}_c^+$ and as such $\J$  we can write.
\[
\J= \left[
\begin{array}{cc}
0 & -G^*\\
G & 0
\end{array}
\right],\;\;GG^*=G^*G={\bf 1}_{{\bS}_c}.
\]
We can view the Hodge-Dolbeault operator  as an operator on ${\bS}_c\oplus {\bS}_c$
\[
\hat{\h}=\J\Biggl(\nah^b_t -\left[
\begin{array}{cc}
\h & 0\\
0 & -G{\h}G^*
\end{array}
\right]\Biggr),\;\; \h^*=\h.
\]
$\h$ is the geometric Dirac operator induced by $\nah^w\otimes {\bf 1}+{\bf 1}\otimes \det \nabla^c$. We want to provide a more explicit description of the operator $\h$. Observe that
\[
C^\infty(\hat{\bS}^+_c)=\Omega^{0,even}(\hat{M})= \Omega^{0,even}(V^*)\oplus \bar{\ve}^0\wedge \Omega^{0,odd}(V^*)
\]
where
\[
\Omega^{0,p}(V^*):= C^\infty(\Lambda^p(V^*)^{0,1}).
\]
We can represent $\psi\in C^\infty(\hat{\bS}^+_c)$ as a sum
\[
\psi=\psi_+\oplus\bar{\ve}^0\wedge \psi_-,\;\;\psi_+\in \Omega^{0,even}(V^*),\;\;\psi_-\in \Omega^{0,odd}(V^*).
\]
The above decomposition   can be alternatively described as follows.   The operator
\[
\bc(\eta)=\J\hbc(\eta):C^\infty(\hat{\bS}^+_c)\ra C^\infty(\hat{\bS}^+_c)
\]
satisfies $\bc(\eta)^2=-1$ and thus $\bc(\ii \eta)$ is an involution of $C^\infty(\hat{\bS}^+_c)$. More explicitly
\[
\bc(\eta)= \frac{\ii}{2}(\hbc(\bar{\ve}^0)+\hbc(\ve^0))(\hbc(\bar{\ve}^0)-\hbc(\ve^0)=\ii(\bar{\ve}^0\wedge\; -\;\ve^0\inpr)(\bar{\ve}^0\wedge\;+ \;\ve^0\inpr).
\]
Thus, for every $\phi\in \Omega^{0,*}(V^*)$ we have
\[
\bc(\ii\eta)(\bar{\ve}^0\wedge \phi)=-\bar{\ve}^0\wedge\phi,\;\;\bc(-\ii \eta)\phi=\phi
\]
This shows that the above decomposition is defined by the $\pm 1$ eigenspaces of the involution $\bc(\eta)$. The restriction of the operator $\bpar: \Omega^{0,*}(\hat{M})\ra \Omega^{0,p}(\hat{M})$ to $\Omega^{0,*}(V^*)$  decomposes into two parts. More precisely, if $\phi\in \Omega^{0,*}(V^*)$ then
\[
\bpar\phi = \bar{\ve}^0\wedge \bpar_0\phi + \bpar_V \phi:=\frac{1}{2}(1+\bc(\ii\eta))\bpar +\frac{1}{2}(1-\bc(\ii\eta))\bpar.
\]
Note that
\[
\bpar_0\phi:= \ve^0\inpr\bpar\phi\in \Omega^{0,p}(V^*),\;\;\bpar_V\in \Omega^{0,p+1}(V^*).
\]
We will regard $\bpar_0$ and $\bpar_V$ as operators
\[
\bpar_0:\Omega^{0,*}(V^*)\ra \Omega^{0,*}(V^*),\;\;\bpar_V:\Omega^{0,*}(V^*)\ra \Omega^{0,*+1}(V^*).
\]
Pick a $t$-independent section $\psi=C^\infty(\hat{\bS}^+_c)$. It decomposes as
\[
\psi=\psi_+ +\bar{\ve}^0\wedge\psi_-,\;\;\psi_\pm \in \Omega^{0,even/odd}(V^*).
\]
We have the equality
\[
\hat{\h}\left[
\begin{array}{c}
\psi\\
0
\end{array}
\right] =- \left[
\begin{array}{cc}
0 & -G^*\\
G & 0
\end{array}
\right]\left[
\begin{array}{cc}
\h & 0\\
0 & -G{\h}G^*
\end{array}
\right]\left[
\begin{array}{c}
\psi\\
0
\end{array}
\right]=\left[
\begin{array}{cc}
0 & \h G^*\\
G\h & 0
\end{array}
\right]\left[
\begin{array}{c}
\psi\\
0
\end{array}
\right]
\]
Thus
\[
\sqrt{2}(\bpar+\bpar^*)\psi =G\h\psi=\hbc(dt)\h \psi\Longrightarrow \h\psi=-\sqrt{2}\J(\bpar+\bpar^*)\psi .
\]
We compute
\[
(\bpar+\bpar^*)(\psi_+ +\bar{\ve}^0\wedge \psi_-) = \bpar\psi_+ +(\bpar\bar{\ve}^0)\wedge \psi_- - \bar{\ve}^0\wedge \bpar\psi_-+\bpar^*\psi_+ +\bpar^*(\bar{\ve}^0\wedge \psi_-)
\]
($\bpar\bar{\ve}^0 =0$)
\[
=\bar{\ve}^0\wedge \bpar_0\psi_+ +\bpar_V\psi_+ -\bar{\ve}^0\wedge \bpar_V\psi_- +(\bar{\ve}^0\wedge \bpar_0 +\bpar_V)^*\psi_+  +\bpar^*(\bar{\ve}^0\wedge \psi_-)
\]
\[
=\bar{\ve}^0\wedge (\bpar_0\psi_+-\bpar_V\psi_-)+\bpar_V\psi_+ +\bpar_V^*\psi_+ +\bpar_0^*(\ve^0\inpr\psi_+)  +\bpar^*(\bar{\ve}^0\wedge \psi_-)
\]
\[
=\bar{\ve}^0\wedge (\bpar_0\psi_+-\bpar_V\psi_-)+\bpar_V\psi_+ +\bpar_V^*\psi_+  +\bpar^*(\bar{\ve}^0\wedge \psi_-).
\]
To proceed further we need to provide a more explicit description  for $\bpar^*(\ve^0\inpr)^*\psi_-$. We  denote by $\lan \bullet,\bullet \ran_M$ the $L^2$-inner product on $M$. For every $t$-independent   compactly supported $\alpha\in \Omega^{0,odd}(\hat{M})$  we have $\alpha=\alpha_- +\bar{\ve}^0\wedge \alpha_+$, $\alpha_\pm \in\Omega^{0,odd/even}(V^*)$, and
\[
\lan \alpha,\bpar^*(\bar{\ve}^0\wedge \phi_-)\ran_M =\lan \bpar\alpha,\bar{\ve}^0\wedge \phi_-\ran_M= \lan\bar{\ve}^0\wedge \bpar_0\alpha_-,\bar{\ve}^0\wedge \phi_-\ran_M -\lan\bar{\ve}^0\wedge \bpar_V\alpha_+,\bar{\ve}^0\wedge\phi_-\ran_M
\]
\[
=\lan\bpar_0\alpha_-, \phi_-\ran_M-\lan \bpar_V\alpha_+,\phi_-\ran_M= \lan \alpha_-,\bpar_0^*\phi_-\ran_M -\lan\alpha_+,\bpar_V^*\phi_-\ran_M
\]
We conclude
\[
\bpar^*(\bar{\ve}^0\wedge \phi_-)= \bpar_0^*\phi_- -\bar{\ve}^0\wedge\bpar_V^*\phi_-,
\]
and
\[
(\bpar+\bpar^*)(\psi_+ +\bar{\ve}^0\wedge \psi_-)=\bar{\ve}^0\wedge (\bpar_0\psi_+-\bpar_V\psi_--\bpar_V^*\phi_-)+\bpar_V\psi_+ +\bpar_V^*\psi_+ +\bpar_0^*\phi_-.
\]
Now observe that
\[
\hbc(dt)\bullet=\frac{1}{\sqrt{2}}(\hbc(\bar{\ve}^0)+\hbc(\ve^0))\bullet=(\bar{\ve}^0\wedge\bullet -\ve^0\inpr\bullet)
\]
so that
\[
\h \psi =-\sqrt{2}(\ve^0\inpr\; -\bar{\ve}^0\wedge)\Bigl\{\bar{\ve}^0\wedge (\bpar_0\psi_+-\bpar_V\psi_--\bpar_V^*\phi_-)+\bpar_V^*\psi_+ +\bpar_V\psi_+ + \bpar_0^*\psi_-\Bigr\}
\]
\[
=-\sqrt{2}\Bigl\{(\bpar_0\psi_+-\bpar_V\psi_--\bpar_V^*\phi_-)-\bar{\ve}^0\wedge(\bpar_V^*\psi_+ +\bpar_V\psi_+ + \bpar_0^*\psi_-)\Bigr\}.
\]
In block form
\[
\h\left[
\begin{array}{c}
\psi_+\\
\psi_-
\end{array}
\right]=\sqrt{2}\left[
\begin{array}{cc}
-\bpar_0 & (\bpar_V^*+\bpar_V)\\
 & \\
 (\bpar_V^*+\bpar_V) & \bpar_0^*
\end{array}
\right]\cdot \left[
\begin{array}{c}
\psi_+\\
\psi_-
\end{array}
\right]
\]
The  above  equality can be further simplified as follows. If $\phi\in \Omega^{0,p}(V^*)\subset \Omega^*(M)\otimes {\bC}$ then
\[
d\phi \in \eta\wedge \Bigl(\Omega^{0,p}(V^*)+\Omega^{1,p-1}(V^*)\Bigr)\oplus \Omega^{0,p+1} (V^*) \oplus \Omega^{1,p}(V^*)\oplus \Omega^{2,p-1}(V^*).
\]
and
\[
-\sqrt{2}\bpar_0 \phi = -\ii (\xi\inpr d\phi)^{0,p} =: -\ii L^V_\xi \phi.
\]
On the other hand, the identity (\ref{eq: bparn}) implies
\[
\bpar_0\phi= \nabla^b_{\bar{\ve}_0}\phi  = \frac{\ii}{\sqrt{2}}\nabla^w_\xi\phi.
\]
Since ${\bf div}_g\xi=0$ the operator $\ii\nabla^w_\xi$ is symmetric  and so must by $\ii L^V_\xi$. Hence $\bpar_0^*\phi=\ii L_\xi^V$ and
\[
\h\left[
\begin{array}{c}
\psi_+\\
\psi_-
\end{array}
\right]= \left[
\begin{array}{cc}
-\ii L^V_\xi  & \sqrt{2}(\bpar_V^*+\bpar_V)\\
 & \\
\sqrt{2}(\bpar_V+\bpar_V^* )& \ii L^V_\xi
\end{array}\right]\cdot \left[
\begin{array}{c}
\psi_+\\
\psi_-
\end{array}
\right]
\]
or equivalently,
\begin{equation}
\h=\bc(\ii \eta)L^V_\xi +\left[
\begin{array}{cc}
0 & \sqrt{2}(\bpar_V+\bpar_V^*)\\
 & \\
\sqrt{2}(\bpar_V^*+\bpar_V^*) & 0
\end{array}\right].
\label{eq: h1}
\end{equation}
We will refer to $\h$ as the {\em  contact Hodge-Dolbeault operator}. The next result summarizes the results we have proved so far.

\begin{theorem} Suppose $(M^{2n+1},g,\eta)$ is a metric contact manifold, $V:=\ker\eta$. Denote by ${\bS}_c$ the bundle of complex spinors  associated to the $spin^c$ structure determined by the contact structure. Denote the corresponding Clifford multiplication by $\bc$.

\medskip

\noindent (i) ${\bS}_c\cong \Lambda^{0,*}V^*$, $\bc(\ii\eta)\phi =(-1)^p\phi$, $\forall \phi\in \Omega^{0,p}(V^*)$.   We decompose
\[
{\bS}_c={\bS}_c^+\oplus {\bS}_c^-,\;\;{\bS}_c^\pm=\Lambda^{0,even/odd}(V^*).
\]
(ii) The operator $\h: C^\infty({\bS}_c)\ra C^\infty({\bS}_c)$   defined by (\ref{eq: h1}) is a geometric Dirac operator induced by the connection $\nabla^w$ on $TM$ and $\det\nabla^c$ on $\det {\bS}_c$.

\medskip

\noindent (iii)  If  we denote by $\dir_c$ the Dirac  operator on ${\bS}_c$ induced by the Levi-Civita connection on $TM$ and $\det \nabla^c$ on $\det {\bS}^c$ then
\[
\h=\dir_c+\frac{1}{4}\bc(\eta\wedge d\eta).
\]
(iv) Using the identity $F(\det\nabla^c)=F(\det\nabla^w)+\frac{n\ii}{2}d\eta$, we deduce that $\h$ satisfies a Weitzenb\"{o}ck formula
\[
\h^2=(\nabla^{\gw})^*(\nabla^{\gw})+\frac{s(g)}{4} +\frac{1}{16}\Bigl(4\bc(d\eta\wedge d\eta)-2n\Bigl) +\frac{1}{2}\bc(F(\det\nabla^w))+\frac{n\ii}{4}\bc(\omega).
\]
In particular, if $\dim M=3$ (so that $n=1$ and $\bc(\eta\wedge d\eta)=-1$) we have
\[
\dir_c=\h+\frac{1}{4},
\]
\[
\h^2=(\nabla^{\gw})^*(\nabla^{\gw})+\frac{s}{4}-\frac{1}{8}+ \frac{1}{2}\bc(F(\det\nabla^w))+\frac{\ii}{4}\bc(d\eta).
\]
\end{theorem}

We want to discuss in more detail the case $\dim M=3$. In this case $\Lambda^{0,even}V^*\cong\uc$ and $\Lambda^{0,odd}(V^*)\cong K_M^{-1}$. The above geometric Dirac operator has the simpler form
\[
\h^2= \bc(\ii \eta)L_\xi +\left[
\begin{array}{cc}
0 & \sqrt{2}\bpar_V^*\\
 & \\
\sqrt{2}\bpar_V & 0
\end{array}\right]=\sqrt{2}\left[
\begin{array}{cc}
 -\bpar_0  & \bpar_V^*\\
 & \\
\bpar_V &  \bpar_0^*
\end{array}
\right]
\]
\[
=\sqrt{2}\left[
\begin{array}{cc}
 -\bpar_0 & 0 \\
 & \\
 0  & \bpar_0^*
 \end{array}
 \right] +\sqrt{2}\left[
 \begin{array}{cc}
 0 & \bpar_V^*\\\\
 & \\
 \bpar_V & 0
 \end{array}
 \right] =: Z+ T.
 \]
Note that along $M$ we have $\bpar_0=\frac{\ii}{\sqrt{2}}\pa_\xi$. We have $\h^2= Z^2+T^2 +\{Z,T\}$, where $\{\bullet,\bullet\}$ denotes the anti-commutator of two operators. In this case
\[
\{Z,T\}= 2\left[
\begin{array}{cc}
0         & [\bpar_0,\bpar_V]^* \\
          &  \\
{[\bpar_0,\bpar_V] }        & 0
\end{array}
\right].
\]
The above commutators can be further simplified using the identity (\ref{eq: bparn}) of \ref{ss: 21}. In this case the Lee $1$-form on $\hat{M}$ is  $dt$.  The equality (\ref{eq: bparn}) implies that for every $t$-independent $\phi\in \Omega^{0,*}(V^*)\subset \Omega^{0,*}(\hat{M})$ we have
\[
\bpar\phi(\bar{\ve}_{k_0},\cdots,\bar{\ve}_{k_p})= \sum_{j=0}^p(-1)^j\nabla^b_{\bar{\ve}_{k_j}}\phi(\bar{\ve}_{k_0},\cdots, \widehat{\bar{\ve}}_{k_j},\cdots \bar{\ve}_{k_p})=\Bigl(\sum_{k=0}^p\bar{\ve}^k\wedge \nabla^b_{\bar{\ve}_k}\Bigr)\phi.
\]
Thus
\[
\bpar_0\phi=\nabla^b_{\bar{\ve}_0}\phi,\;\;\bpar_V \phi =\Bigl(\sum_{k=1}^p\bar{\ve}^k\wedge \nabla^b_{\bar{\ve}_k}\Bigr)\phi.
\]
When $\dim M=3$    and   $\phi=u\in \Omega^{0,0}(V^*)=C^\infty(M)\otimes {\bC}$ we have
\[
[\bpar_0,\bpar_V]u=\nabla^b_{\bar{\ve}_0}(\bar{\ve}^1\wedge \nabla^b_{\bar{\ve}_1} u)-\bar{\ve}^1\wedge \nabla^b_{\bar{\ve}_1}\nabla^b_{\bar{\ve}_0}u
\]
\[
= (\nabla^b_{\bar{\ve}_0}\bar{\ve}^1)\wedge \nabla^b_{\bar{\ve}_1}
u+\bar{\ve}^1\wedge [\nabla^b_{\bar{\ve}_0}, \nabla^b_{\bar{\ve}_1}]u=(\nabla^b_{\bar{\ve}_0}\bar{\ve}^1)\wedge \nabla^b_{\bar{\ve}_1} u +\bar{\ve}^1\wedge\nabla^b_{[\bar{\ve}_0,\bar{\ve}_1]}u +\bar{\ve}_1 \wedge F_b(\bar{\ve}_0,\bar{\ve}_1)u,
\]
where $F_b$ denotes the curvature of the $\nabla^b$. Denote by $T_b$ the torsion of $\nabla^b$. Observe that
\[
\nabla^b_{\bar{\ve}_0}\bar{\ve}^1(\bar{\ve}_1)=-\bar{\ve}^1(\nabla^b_{\bar{\ve}_0}\bar{\ve}_1)
\]
so that
\[
(\nabla^b_{\bar{\ve}_0}\bar{\ve}^1)\wedge \nabla^b_{\bar{\ve}_1} u+\bar{\ve}^1\wedge\nabla^b_{[\bar{\ve}_0,\bar{\ve}_1]}u=-\bar{\ve}^1\wedge \nabla^b_{\nabla^b_{\bar{\ve}_0}\bar{\ve}_1}u+\bar{\ve}^1\wedge\nabla^b_{[\bar{\ve}_0,\bar{\ve}_1]}u
\]
($\nabla^b\bar{\ve}_0=0$)
\[
=-\bar{\ve}^1\wedge(\nabla^b_{\nabla^b_{\bar{\ve}_0}\bar{\ve}_1}-\nabla^b_{\nabla^b_{\bar{\ve}_1}\bar{\ve}_0}-\nabla^b_{[\bar{\ve}_0,\bar{\ve}_1]})u=\bar{\ve}^1\wedge \nabla^b_{T_w(\bar{\ve}_1,\bar{\ve}_0)}u=\frac{\ii}{\sqrt{2}}\bar{\ve}^1\wedge\nabla^w_{T_w(\bar{\ve}_1,\xi)}u
\]
\[
=-\frac{\ii}{2\sqrt{2}}\bar{\ve}^1\wedge \nabla^b_{J\Phi\bar{\ve}_1}u=\frac{\ii}{2\sqrt{2}}\bar{\ve}^1\wedge \nabla^b_{\Phi J\bar{\ve}_1}u
\]
($J\bar{\ve}-1=-\ii\bar{\ve}_1$)
\[
= \frac{1}{2\sqrt{2}}\bar{\ve}^1\wedge \nabla^b_{\Phi\bar{\ve}_1}u =\frac{1}{2\sqrt{2}}\Phi_c\bigl(\pa_Vu \bigr)=:{\gT}u.
\]
where $\Phi_c$ is  the  complexification\footnote{$\Phi_c$ is complex linear  but it {\em anticommutes} with $J$.} of $\Phi$.  The differential operator by ${\gT}$ is trivial when $\Phi=0$, which in the $3$-dimensional case is equivalent to $M$ being  Sasakian or to $\hat{J}$ being integrable.

Putting together all the above facts we obtain
\[
[\bpar_0, \bpar_V]=\bar{\ve}^1\wedge F_b(\bar{\ve}_0,\bar{\ve}_1) +{\gT}.
\]
We conclude
\[
\{Z, T\} =2\left[
\begin{array}{cc}
0 & \overline{F_b(\bar{\ve}_0,\bar{\ve}_1)}\ve^1\inpr\\
& \\
\bar{\ve}^1\wedge F_b(\bar{\ve}_0,\bar{\ve}_1) & 0
\end{array}
\right]+ 2\tilde{\gT},\;\;\tilde{\gT}:=\left[
\begin{array}{cc}
0  & {\gT}^*\\
 & \\
{\gT}& 0
 \end{array}
 \right].
 \]
The zero order operator  above can be further simplified by observing that
\[
\bc(\bar{\ve}^1)=\sqrt{2}\left[
\begin{array}{cc}
0 & 0 \\
\bar{\ve}^1\wedge & 0
\end{array}
\right],\;\; \bc(\ve^1)=\sqrt{2}\left[
\begin{array}{cc}
0 & -\ve^1\inpr \\
0 & 0
\end{array}
\right]
\]
so that
\[
\left[
\begin{array}{cc}
0 & \overline{F_b(\bar{\ve}_0,\bar{\ve}_1)}\ve^1\inpr\\
& \\
\bar{\ve}^1\wedge F_b(\bar{\ve}_0,\bar{\ve}_1) & 0
\end{array}
\right]=\frac{1}{\sqrt{2}}\bc\Bigl(F_b(\bar{\ve}_0,\bar{\ve}^1)\bar{\ve}^1-\overline{F_b(\bar{\ve}_0,\bar{\ve}_1)}\ve^1\Bigr)
\]
\[
=\frac{\ii}{2}\bc\Bigl(F_w(\xi,\bar{\ve}_1)\bar{\ve}^1+\overline{F_w(\xi,\bar{\ve}_1)}\ve^1\Bigr)
\]
Above we denoted by $F_w$ the curvature of $\nabla^w$ as a connection  on the hermitian line bundle $(V,J)\cong K_M^{-1}$. To get a more suggestive description  we write
\[
\xi\inpr F_w= \ii(ae^1 +bf^1),
\]
where $a,b$ are locally defined {\em real valued functions}. Then
\[
F_w(\xi,\bar{\ve}_1)\bar{\ve}^1= \frac{\ii}{2}(a+\ii b) (e^1-\ii f^1),\;\;\overline{F_w(\xi,\bar{\ve}_1)}\ve^1=-\frac{\ii}{2}(a-\ii b)(e^1+\ii f^1).
\]
Thus
\[
F_w(\xi,\bar{\ve}_1)\bar{\ve}^1+\overline{F_w(\xi,\bar{\ve}_1)}\ve^1=(-be^1+af^1)=-\ii (\ast F_w -\eta\wedge(\xi\inpr \ast F_w)).
\]
The last term can also be described as $-\ii P_V (\ast F_w)$, where  $P_V$ denotes the orthogonal projection $T^M\ra V^*$, and  $\ast$  denotes the {\em complex linear} extension of the Hodge operator. The above  facts   now yield the following commutator identities.
\begin{subequations}
\begin{equation}
\{ Z, T\}= \bc\bigl(P_V \ast F_w\bigr) + 2\tilde{\gT},
\label{eq: ant-comm}
\end{equation}
\begin{equation}
\h^2=Z^2+T^2 + \bc\bigl(P_V \ast F_w\bigr) + 2\tilde{\gT}.
\label{eq: decouple}
\end{equation}
\label{subeq: com}
\end{subequations}

\begin{remark} (a) If we twist the Dirac operator $\dir(\nabla^w)$ by  a  hermitian connection on the trivial line bundle $\uc$ we obtain a  new Dirac operator $\h_A$ satisfying
\[
\h_A=\left[
\begin{array}{cc}
-\ii \nabla^A_\xi  & \sqrt{2}(\bpar^A_V)^*\\
 & \\
\sqrt{2}\bpar^A_V& \ii \nabla^A_\xi
\end{array}\right]=:Z_A+T_A.
\]
The operators $Z_A$ and $T_A$ satisfy the anticommutation rule
\begin{equation}
\{Z_A,T_A\}=Z_A^2+ T_A^2 + \bc\bigl(P_V\ast F_w)\bigr) + \bc\bigl(P_V\ast  F_A\bigr)+ 2\tilde{\gT}_A
\label{eq: ant-coup}
\end{equation}
where $\tilde{\gT}_A$ is defined  as $\tilde{\gT}$ using instead the operator ${\gT}_A :=\frac{1}{2\sqrt{2}}\Phi_c\pa^A_V$.

\medskip

\noindent (b) The curvature $F_w$ has the local description
\[
F_w  =-\ii \rho d\eta + \eta\wedge (\xi\inpr F_w).
\]
Up to a positive multiplicative constant (depending on various normalization conventions) the scalar   $\rho$ is known as the  {\em Webster scalar curvature}. We refer to \cite{CH} for more details.
\end{remark}

\subsection{Connections induced by symplectizations}
\label{ss: 33}

The  symplectization of the   positively oriented metric contact
manifold $(M^{2n+1},\eta, g, J)$ is the manifold $\tilde{M}= {\bR}_+\times M$ equipped with the symplectic form
\[
\tilde{\omega}=dt\wedge \eta + t d\eta=dt\wedge \eta +t\omega.
\]
If we denote by $\td$ the exterior derivative on $\tilde{M}$ then we can write
\[
\tilde{\omega}= \td(t\eta).
\]
$\tilde{M}$ is equipped with a compatible  metric
\[
\tilde{g}= dt^2 +\eta^2 +t\omega(\bullet,J\bullet).
\]
We denote by $\tilde{J}$ the associated almost complex structure. We will identify  $M$ with the slice $\{1\} \times M$ of $\tilde{M}$.

If we fix as before a local, oriented, orthonormal frame $\xi,e_1,f_1,\cdots, e_n ,f_n$ compatible with the metric contact structure on $M$ then we get a symplectic frame
\[
\te_0=\pa_t,\;\tf_0=\xi,\;\;\te_k=t^{-1/2}e_k,\;\;\tf_k= t^{-1/2}f_k,\;\;k=1,\cdots, n.
\]
The dual coframe is
\[
\te^0 =dt, \; \tf^0=\eta,\;\; \te^k=t^{1/2}e^k,\;\;\tf^k=t^{1/2}f^k, \;\;k=1,\cdots, n.
\]
We denote by $\tilde{N}$ the Nijenhuis tensor of $\tilde{J}$ and by $\hN$ the Nijenhuis tensor of the  almost complex manifold $(\hat{M},\hat{J})$ used in \ref{ss: 31}.  The    Chern connection  $\tna^c$ of $(\tilde{M}, \tilde{g},\tilde{J})$ is the  metric connection with torsion $\tilde{T} =\tilde{N}$. In this case
\[
\theta=0,\;\;{\gb}\tilde{T}=0.
\]
Observe that  $\tilde{J}=\hat{J}$. We deduce   that for $j,k=1,\cdots ,n$ we have
\[
\tilde{N}(\te_j,\te_k)= \frac{1}{t}\hN(e_j,e_k),\;\; \tilde{N}(\te_j, \tf_k)=\frac{1}{t}\hN(e_j,f_k),\;\;\tilde{N}(\tf_j,\tf_k)=\frac{1}{t}\hN(f_j,f_k),
\]
\[
\tilde{N}(\pa_t, \te_j)=\frac{1}{\st}\hN(\pa_t, e_j),\;\; \tilde{N}(\pa_t,\tf_k)=\frac{1}{\st} \hN(\pa_t, f_k),
\]
\[
\tilde{N}(\xi, \te_j)=\frac{1}{\st}\hN(\pa_t, e_j),\;\; \tilde{N}(\xi_t,\tf_k)=\frac{1}{\st} \hN(\pa_t, f_k).
\]
Denote by $\tilde{D}$ the  Levi-Civita connection determined by $\tg$. It determined by (see \cite{N1})
\[
2\tg (\tilde{D}_X Y, Z)=X\tg(Y,Z)+Y\tg(X,Z)-Z\tg(X,Y)
\]
\[
 +\tg([X,Y],Z)+\tg([Z,X],Y)+\tg(X,[Z,Y]).
\]
We deduce from the above identity that if $X,Y$ are $t$-independent  vectors tangent along $M$
\[
2\tg(\tilde{D}_tX,Y)=g(X_V,Y_V)=\omega(X,JY),
\]
where $X_V:= P_VX$.
\[
2\tg(\tilde{D}_XY,\pa_t)= -\pa_t\tg(X,Y)= -g(X_V,Y_V) = \omega(JX,Y).
\]
As in \ref{ss: 31} we want to alter $\tna^c$ by  $B\in \Omega^{1,1}_s(T^*\tilde{M})$ such that $\tr B=0$ so that the new  basic hermitian connection  $\tna^b$ with torsion $\tilde{T}_b^\dag:=\tilde{N}^\dag +B$ satisfies
\begin{equation}
\tna^b_X\xi=0, \;\;\tg(\tna^b_XY, \pa_t)=0,
\label{eq: condit}
\end{equation}
for all  $t$-independent tangent vectors $X$, $Y$ along $M$.

We have $\tna= \tilde{D}+A$, where $A^\dag =-\tilde{T}_b^\dag$. Thus we  need
\[
0=\tg(\tna_X Y, \pa_t)=\tg(\tilde{D}_XY,\pa_t))- \tg(X, \tilde{N}(Y, \pa_t)-B(X; Y,\pa_t)
\]
\[
=-\frac{1}{2}\omega(JX,Y)+\tg(X, \tilde{N}(\pa_t, Y))- B(X; Y,\pa_t)
\]
If $Y=\xi$ we deduce
\[
B(X; \xi,\pa_t)=0.
\]
If $Y\in C^\infty(V)$ then we deduce
\[
0=-\frac{1}{2}\omega(JX,Y)+\frac{1}{\st}\tg(X,\hN(\pa_t, Y))+B(X;\pa_t,Y)
\]
\[
=\frac{1}{2}g(X,Y)+\frac{1}{4\st}\tg(X,\Phi Y) +B(X;\pa_t,Y)=\frac{1}{2}g(X,Y)+\frac{\st}{4}g(X,\Phi Y)+B(X;\pa_t,Y)
\]
We conclude that $B$ must satisfy the additional conditions
\[
B(\xi;\pa_t, Y)=0,\;\;Y\in C^\infty(V)
\]
\[
B(X;\pa_t, Y)=-\frac{1}{2\st}\Bigl( \frac{1}{\st}\tg(X, Y)+\frac{1}{2}\tg(X,\Phi Y)\Bigr)
\]
We write $B= B_0+ B_1$ where $B_0$ is defined  as in Lemma \ref{lemma: b0} by the equality
\[
B_0= \frac{1}{4\st}\Bigl\{ \Phi\wedge dt+  J\Phi\wedge \eta\Bigr\}
\]
$B_1$ must satisfy the equalities $\tr B_1=0$,
\begin{subequations}
\begin{equation}
B_1(X;\pa_t, Y)=-\frac{1}{2t}\tg(X,Y),\;\;\forall X,Y\in C^\infty(V)
\label{eq: b1a}
\end{equation}
\begin{equation}
B_1(X; \xi,\pa_t)=B_1(\xi;\pa_t, Y)=0,\;\;\forall X\in \Vect(M),\;Y\in C^\infty(V)
\label{eq: b1b}
\end{equation}
\end{subequations}
We try $B_1$ of the form
\[
B_1= xdt\otimes dt\wedge \eta + y \eta \otimes d\eta +  U +V
\]
where
\[
U= \frac{1}{2t}P_V\wedge dt ,\;\;V=\frac{1}{2t}JP_v\wedge \eta
\]
Clearly  $B_1\in \Omega^{1,1}(T^*\tilde{M})$. Next observe that
\[
{\gb}B_1= y\eta\wedge d\eta + {\gb}V= (y +\frac{1}{t})\eta\wedge d\eta,
\]
\[
\tr B_1 = (x +\frac{n}{t})dt,
\]
Thus, set $x=-\frac{n}{t}$, $y=\frac{1}{t}$. These choices guarantee that $B_1\in \Omega^{1,1}_s(T^*\tilde{M})$ and $\tr B_1= 0$. The conditions  (\ref{eq: b1a})  and (\ref{eq: b1b}) can now be verified  by direct computation.     We can now conclude that if
\[
B=\frac{1}{4\st}\Bigl(\Phi\wedge dt+J\Phi\wedge \eta\Bigr) -\frac{n}{t}dt\otimes dt\wedge \eta-\frac{1}{t} \eta \otimes d\eta + \frac{1}{2t}\Bigl(P_V\wedge dt +JP_V\wedge \eta\Bigr)
\]
then the connection $\tna^b$ with torsion $\tilde{N}^\dag +B$ satisfies the conditions  (\ref{eq: condit}). These conditions show that $\tna^b$ induces,,  by restriction to each slice $\{t\}\times M$, a  connection $\nabla^t$ on $TM$.     The torsion of $\nabla^1=\nabla^{t=1}$ is given by
\[
(T_1)^\dag = \tilde{N}^\dag\mid_{t=1}+ B\mid_{t=1}=  \hN^\dag\mid_{M}+\frac{1}{4}(J\Phi\wedge \eta)-\eta\otimes d\eta +\frac{1}{2}(JP_V\wedge \eta)
\]
\[
=\frac{1}{2}N^\dag -\frac{1}{2}\eta\otimes d\eta +\frac{1}{2}(JP_V\wedge \eta)+\frac{1}{4}(J\Phi\wedge \eta).
\]
When $M$ is a $CR$ manifold we deduce
\[
T_1^\dag=-\eta\otimes d\eta +\frac{1}{2}J\wedge \eta+\frac{1}{2}J\Phi\wedge \eta.
\]
In particular
\[
T^1(X,Y)=-\xi d\eta(X,Y),\;\;\forall X,Y\in C^\infty(V).
\]
This connection  {\em never}  coincides with generalized Webster connection constructed in \ref{ss: 31}, because in this case we have ${\gb}T_1^\dag =0$. This shows  $\nabla^1$ is Dirac equivalent to  the Levi-Civita connection.  We have thus proved  the following result.

\begin{theorem} On every    metric contact manifold  $(M,g,J)$ there exists   a  canonical   nice  contact metric connection $=\nabla^1$ induced  by a  basic Hermitian connection on the symplectization  of $M$. This contact connection is  Dirac equivalent to the Levi-Civita connection and its torsion is given by
\[
T_1^\dag= \frac{1}{2}N^\dag -\frac{1}{2}\eta\otimes d\eta +\frac{1}{2}(JP_V\wedge \eta)+\frac{1}{4}(J\Phi\wedge \eta).\;\;\Box
\]
\label{th: sympl}
\end{theorem}

Let us observe that if $M$ is $CR$ then for every $X,Y\in C^\infty(V)$ we have
\[
g(\nabla^1_\xi X, Y)= -g(\xi, T(X,Y))= +\omega(X,Y)
\]
so that
\[
\nabla^1_\xi  = D^V_\xi-J= P_VD_\xi+J=\nabla^w_\xi+J.
\]

\begin{proposition} Suppose $M$ is a $CR$ manifold. Then
\[
\det \nabla^1 =\det \nabla^w +3n\ii\eta=\det\nabla^c+\frac{5n\ii}{2}\eta.
\]
\label{prop: compare1}
\end{proposition}

\noindent {\bf Proof}\hspace{.3cm}
\[
\Delta:=T_1^\dag-T_w^\dag =\frac{1}{4}(-7\eta\otimes d\eta -\eta\wedge d\eta +3J\wedge \eta)
\]
so that $\nabla^1=\nabla^w+ A$ where
\[
A^\dag= \frac{1}{2}{\gb}\Delta -\Delta =-\frac{1}{2}\eta\wedge d\eta+ \frac{1}{4}(7\eta\otimes d\eta +\eta\wedge d\eta -3J\wedge \eta)
\]
\[
= \frac{1}{4}(7\eta\otimes d\eta -\eta\wedge d\eta -3J\wedge \eta).
\]
Set
\[
\delta=\ve_1\wedge \cdots \wedge \ve_n.
\]
As in the proof  of  Proposition \ref{prop: compare} we have
\[
\det \nabla^1_X \delta =\det \nabla^w_X  +\ii\bigl(\sum_{k=1}^nC_k^k\bigr)\delta
\]
where
\[
C_k^k= g(A_Xe_k, Jf_k) +g(A_X f_k, Jf_k).
\]
The above sum is nontrivial only for $X=\xi$ in which case it is equal to $3n$. We conclude that
\[
\det \nabla^1=\det \nabla^w + 3n\ii \eta.\;\;\Box
\]

\begin{remark} Let us point out  a difference between contact and Hermitian connections.  We have shown that there always exist contact connections   with torsion $T$ satisfying ${\gb}T^\dag =0$.

On the other hand, if   $\nabla$ is a Hermitian connection  on an almost complex Hermitian manifold $(M,g,J)$ with Nijenhuis tensor $N$  then  its torsion satisfies (see \cite{Gau})
\[
({\gb}T)^-=({\gb}N^\dag)=(d^c\omega)^-.
\]
If $\dim M =4$ then  always $(d^c\omega)^-=0$ and  in this case it is possible to find Hermitian connections Dirac equivalent to the Levi-Civita connection.  However, in higher dimensions this is possible if and only if $(d^c\omega)^-=0$.
\end{remark}

\subsection{Uniqueness results}
\label{ss: 34}

The constructions we  performed in the previous subsection may seem  a bit ad-hoc   but as we will show   in  this section they produce, at least  for $CR$ manifolds, connections uniquely determined by a few natural requirements.

\begin{proposition} Suppose $(M, \eta,g, J)$ is  $CR$ connection. Then  each  Dirac equivalence class of connections contains at most one nice $CR$ connection.
\label{prop: uniq}
\end{proposition}

\noindent {\bf Proof}\hspace{.3cm}  Suppose $\nabla$ is a nice  $CR$ connection with torsion $T$. Set $\Omega:={\gb}T$. We get a hermitian connection $\nah=dt\wedge\pa_t+\nabla$  on $(T\hat{M}, \hat{J})$ with the property
\[
{\gb}T(\nah)^\dag = \Omega,\;\; \tr T(\nah)^\dag=0.
\]
Denote by $\nabla^b$ the  basic hermitian connection on $(T\hat{M},\hat{J})$ we have constructed in  \ref{ss: 31}.   The results in \ref{ss: 21} imply that
\[
T(\nah)^\dag =T_b^\dag + \frac{9}{8} \psi^+-\frac{3}{8}{\gM}\psi^+ + B =: T_b^\dag + S,
\]
where
\[
\psi^+\in \Omega^{3,+}(\hat{M}),\;\;B\in \Omega^{1,1}_s(T^*\hat{M}),
\]
\[
\Omega={\gb}T_b^\dag +3\psi^+=3\psi^+ +\eta\wedge d\eta,
\]
\begin{equation*}
B(\pa_t;\bullet,\bullet)=0= B(\bullet;\bullet,\pa_t)=0,\;\;\tr B=0.
\tag{$\ast$}
\label{tag: ast}
\end{equation*}
Thus $\psi^+$ is uniquely determined. Moreover, since $\nabla$ is a $CR$ connection   we deduce  that
\[
g(X,T(Y,Z))=0,\;\; \forall X,Y,Z\in C^\infty(V).
\]
Since the restriction of $\nabla^b$ to $M$ is also a $CR$ connection we deduce
\[
S(X;Y,Z)=0,\;\;\forall X,Y,Z\in C^\infty(V).
\]
Thus the restriction of $B$ to $V$ is uniquely determined.     The  condition $B\in \Omega^{1,1}_s(T^*\hat{M})$  coupled with (\ref{tag: ast}) show that the restriction of $B$ to ${\bR}\pa_t\oplus {\bR}\xi\subset T\hat{M}$ is  also uniquely determined. This concludes the proof of Proposition \ref{prop: uniq}.$\Box$

\bigskip

\begin{remark} We can use  Gauduchon's description   of   the  hermitian connections on $T\hat{M}$ to completely characterize which   Dirac equivalence classes of connections on $TM$ contain nice  $CR$ connections.
\end{remark}

\begin{corollary} The Webster connection on a $CR$ manifold is the unique  $CR$ connection adapted to $\h$.  Moreover, the connection $\nabla^1$ of  \ref{ss: 33} is the unique nice $CR$ connection    with  torsion satisfying ${\gb}T^\dag=0$.
\end{corollary}

\end{document}